\documentclass[10pt]{amsart}
\usepackage{amsmath}
\usepackage{amsxtra}
\usepackage{amscd}
\usepackage{amsthm}
\usepackage{amsfonts}
\usepackage{amssymb}
\usepackage{eucal}

\newtheorem{corollary}[subsection]{Corollary}
\newtheorem{lemma}[subsection]{Lemma}
\newtheorem{proposition}[subsection]{Proposition}

\newtheorem{theorem}[subsection]{Theorem}
\newtheorem{definition}[subsection]{Definition}
\newtheorem{remark}[subsection]{Remark}

\newcommand{\nc}{\newcommand}
\nc{\renc}{\renewcommand} \nc{\ssec}{\subsection}
\nc{\sssec}{\subsubsection} \nc{\on}{\operatorname}

\nc\ol{\overline} \nc\ul{\underline} \nc\wt{\widetilde}
\nc\tboxtimes{\wt{\boxtimes}} \nc{\alp}{\alpha}

\nc{\ZZ}{{\mathbb Z}} \nc{\NN}{{\mathbb N}} \nc{\CC}{{\mathbb C}}
\nc{\OO}{{\mathbb O}} \renc{\SS}{{\mathbb S}} \nc{\DD}{{\mathbb D}}

\nc{\Fq}{{\mathbb F}_q} \nc{\Fqb}{\ol{{\mathbb F}_q}}
\nc{\Ql}{\ol{{\mathbb Q}_\ell}} \nc{\id}{\text{id}} \nc\X{\mathcal
X}

\nc{\Hom}{\on{Hom}} \nc{\Lie}{\on{Lie}} \nc{\Loc}{\on{Loc}}
\nc{\Pic}{\on{Pic}} \nc{\Bun}{\on{Bun}} \nc{\IC}{\on{IC}}
\nc{\Aut}{\on{Aut}} \nc{\rk}{\on{rk}} \nc{\Sh}{\on{Sh}}
\nc{\Perv}{\on{Perv}} \nc{\pos}{{\on{pos}}} \nc{\Conv}{\on{Conv}}
\nc{\Sph}{\on{Sph}} \nc{\Sym}{\on{Sym}}
\nc{\BunBb}{\overline{\Bun}_B} \nc{\Buno}{\overset{o}{\Bun}}
\nc{\BunPb}{{\overline{\Bun}_P}} \nc{\BunBM}{\overline{\Bun}_{B(M)}}
\nc{\BunPbw}{{\widetilde{\Bun}_P}}
\nc{\BunBP}{\widetilde{\Bun}_{B,P}} \nc{\GUb}{\overline{G/U}}
\nc{\GUPb}{\overline{G/U(P)}}
\nc{\iso}{{\stackrel{\sim}{\longrightarrow}}}

\nc{\Hhom}{\underline{\on{Hom}}} \nc\syminfty{\on{Sym}^{\infty}}
\nc\lal{\ol{\lambda}} \nc\xl{\ol{x}} \nc\thl{\ol{\theta}}
\nc\nul{\ol{\nu}} \nc\mul{\ol{\mu}} \nc\Sum\Sigma
\nc{\oX}{\overset{o}{X}{}}

\nc{\M}{{\mathcal M}} \nc{\N}{{\mathcal N}} \nc{\F}{{\mathcal F}}
\nc{\D}{{\mathcal D}} \nc{\Q}{{\mathcal Q}} \nc{\Y}{{\mathcal Y}}
\nc{\G}{{\mathcal G}} \nc{\E}{{\mathcal E}} \nc{\CalC}{{\mathcal C}}
\nc\Dh{\widehat{\D}}

\nc{\C}{{\mathcal C}} \nc{\K}{{\mathcal K}}
\renewcommand{\H}{{\mathcal H}}

\nc{\T}{{\mathcal T}} \nc{\V}{{\mathcal V}} \renc{\P}{{\mathcal P}}
\nc{\A}{{\mathcal A}} \nc{\B}{{\mathcal B}} \nc{\U}{{\mathcal U}}

\nc{\Gr}{\on{Gr}}

\nc{\frn}{{\check{\mathfrak u}(P)}}
\nc\f{{\mathfrak f}}

\nc\w{\text{w}}

\nc\mathi\iota \nc\Spec{\on{Spec}} \nc\Mod{\on{Mod}}
\nc{\tw}{\widetilde{\mathfrak t}} \nc{\pw}{\widetilde{\mathfrak p}}
\nc{\qw}{\widetilde{\mathfrak q}} \nc{\jw}{\widetilde j}

\nc{\grb}{\overline{\Gr}} \nc{\I}{\mathcal I}

\nc{\lambdach}{{\check\lambda}} \nc{\Lambdach}{{\check\Lambda}{}}
\nc{\much}{{\check\mu}} \nc{\omegach}{{\check\omega}}
\nc{\nuch}{{\check\nu}} \nc{\etach}{{\check\eta}}
\nc{\alphach}{{\check\alpha}} \nc{\betach}{{\check\beta}}
\nc{\rhoch}{{\check\rho}} \nc{\ch}{{\check h}}

\nc{\Hb}{\overline{\H}}

\nc{\BA}{{\mathbb{A}}} \nc{\BC}{{\mathbb{C}}} \nc{\BM}{{\mathbb{M}}}
\nc{\BN}{{\mathbb{N}}} \nc{\BP}{{\mathbb{P}}} \nc{\BR}{{\mathbb{R}}}
\nc{\BZ}{{\mathbb{Z}}} \nc{\BS}{{\mathbb{S}}}

\nc{\CA}{{\mathcal{A}}} \nc{\CB}{{\mathcal{B}}}
\nc{\CE}{{\mathcal{E}}} \nc{\CF}{{\mathcal{F}}}
\nc{\CG}{{\mathcal{G}}} \nc{\CH}{{\mathcal{H}}}
\nc{\CI}{{\mathcal{I}}} \nc{\CL}{{\mathcal{L}}}
\nc{\CM}{{\mathcal{M}}} \nc{\CN}{{\mathcal{N}}}
\nc{\CO}{{\mathcal{O}}} \nc{\CP}{{\mathcal{P}}}
\nc{\CQ}{{\mathcal{Q}}} \nc{\CR}{{\mathcal{R}}}
\nc{\CS}{{\mathcal{S}}} \nc{\CT}{{\mathcal{T}}}
\nc{\CU}{{\mathcal{U}}} \nc{\CV}{{\mathcal{V}}}
\nc{\CW}{{\mathcal{W}}} \nc{\CZ}{{\mathcal{Z}}}

\nc{\cM}{{\check{\mathcal M}}{}} \nc{\csM}{{\check{\mathcal A}}{}}
\nc{\oM}{{\overset{\circ}{\mathcal M}}{}}
\nc{\obM}{{\overset{\circ}{\mathbf M}}{}}
\nc{\oCA}{{\overset{\circ}{\mathcal A}}{}}
\nc{\obA}{{\overset{\circ}{\mathbf A}}{}}
\nc{\ooM}{{\overset{\circ}{M}}{}} \nc{\osM}{{\overset{\circ}{\mathsf
M}}{}} \nc{\vM}{{\overset{\bullet}{\mathcal M}}{}}
\nc{\nM}{{\underset{\bullet}{\mathcal M}}{}}
\nc{\oD}{{\overset{\circ}{\mathcal D}}{}}
\nc{\obD}{{\overset{\circ}{\mathbf D}}{}}
\nc{\oA}{{\overset{\circ}{\mathbb A}}{}}
\nc{\op}{{\overset{\bullet}{\mathbf p}}{}}
\nc{\cp}{{\overset{\circ}{\mathbf p}}{}}
\nc{\oU}{{\overset{\bullet}{\mathcal U}}{}}
\nc{\ofZ}{{\overset{\circ}{\mathfrak Z}}{}}

\nc{\ff}{{\mathfrak{f}}} \nc{\fv}{{\mathfrak{v}}}
\nc{\fa}{{\mathfrak{a}}} \nc{\fb}{{\mathfrak{b}}}
\nc{\fd}{{\mathfrak{d}}} \nc{\fe}{{\mathfrak{e}}}
\nc{\fg}{{\mathfrak{g}}} \nc{\fgl}{{\mathfrak{gl}}}
\nc{\fh}{{\mathfrak{h}}} \nc{\fri}{{\mathfrak{i}}}
\nc{\fj}{{\mathfrak{j}}} \nc{\fk}{{\mathfrak{k}}}
\nc{\fm}{{\mathfrak{m}}} \nc{\fn}{{\mathfrak{n}}}
\nc{\ft}{{\mathfrak{t}}} \nc{\fu}{{\mathfrak{u}}}
\nc{\fw}{{\mathfrak{w}}} \nc{\fz}{{\mathfrak{z}}}
\nc{\fp}{{\mathfrak{p}}} \nc{\frr}{{\mathfrak{r}}}
\nc{\fs}{{\mathfrak{s}}} \nc{\fsl}{{\mathfrak{sl}}}
\nc{\hsl}{{\widehat{\mathfrak{sl}}}}
\nc{\hgl}{{\widehat{\mathfrak{gl}}}}
\nc{\hg}{{\widehat{\mathfrak{g}}}}
\nc{\chg}{{\widehat{\mathfrak{g}}}{}^\vee}
\nc{\hn}{{\widehat{\mathfrak{n}}}}
\nc{\chn}{{\widehat{\mathfrak{n}}}{}^\vee}

\nc{\fA}{{\mathfrak{A}}} \nc{\fB}{{\mathfrak{B}}}
\nc{\fD}{{\mathfrak{D}}} \nc{\fE}{{\mathfrak{E}}}
\nc{\fF}{{\mathfrak{F}}} \nc{\fG}{{\mathfrak{G}}}
\nc{\fH}{{\mathfrak{H}}} \nc{\fI}{{\mathfrak{I}}}
\nc{\fJ}{{\mathfrak{J}}} \nc{\fK}{{\mathfrak{K}}}
\nc{\fL}{{\mathfrak{L}}} \nc{\fM}{{\mathfrak{M}}}
\nc{\fN}{{\mathfrak{N}}} \nc{\frP}{{\mathfrak{P}}}
\nc{\fQ}{{\mathfrak{Q}}} \nc{\fT}{{\mathfrak{T}}}
\nc{\fU}{{\mathfrak{U}}} \nc{\fV}{{\mathfrak{V}}}
\nc{\fW}{{\mathfrak{W}}} \nc{\fX}{{\mathfrak{X}}}
\nc{\fY}{{\mathfrak{Y}}} \nc{\fZ}{{\mathfrak{Z}}}

\nc{\ba}{{\mathbf{a}}} \nc{\bb}{{\mathbf{b}}} \nc{\bc}{{\mathbf{c}}}
\nc{\be}{{\mathbf{e}}}
\nc{\bff}{{\mathbf{f}}}\nc{\bk}{{\mathbf{k}}}\nc{\bj}{{\mathbf{j}}}
\nc{\bn}{{\mathbf{n}}} \nc{\bp}{{\mathbf{p}}} \nc{\bq}{{\mathbf{q}}}
\nc{\br}{{\mathbf{r}}} \nc{\bfu}{{\mathbf{u}}}
\nc{\bv}{{\mathbf{v}}} \nc{\bx}{{\mathbf{x}}} \nc{\by}{{\mathbf{y}}}
\nc{\bw}{{\mathbf{w}}} \nc{\bA}{{\mathbf{A}}} \nc{\bB}{{\mathbf{B}}}
\nc{\bC}{{\mathbf{C}}} \nc{\bD}{{\mathbf{D}}} \nc{\bF}{{\mathbf{F}}}
\nc{\bH}{{\mathbf{H}}} \nc{\bK}{{\mathbf{K}}} \nc{\bM}{{\mathbf{M}}}
\nc{\bN}{{\mathbf{N}}} \nc{\bO}{{\mathbf{O}}} \nc{\bS}{{\mathbf{S}}}
\nc{\bV}{{\mathbf{V}}} \nc{\bW}{{\mathbf{W}}} \nc{\bX}{{\mathbf{X}}}
\nc{\bY}{{\mathbf{Y}}} \nc{\bP}{{\mathbf{P}}} \nc{\bZ}{{\mathbf{Z}}}
\nc{\bh}{{\mathbf{h}}}

\nc{\sA}{{\mathsf{A}}} \nc{\sB}{{\mathsf{B}}} \nc{\sC}{{\mathsf{C}}}
\nc{\sD}{{\mathsf{D}}} \nc{\sE}{{\mathsf{E}}} \nc{\sF}{{\mathsf{F}}}
\nc{\sK}{{\mathsf{K}}} \nc{\sL}{{\mathsf{L}}} \nc{\sM}{{\mathsf{M}}}
\nc{\sO}{{\mathsf{O}}} \nc{\sQ}{{\mathsf{Q}}} \nc{\sP}{{\mathsf{P}}}
\nc{\sT}{{\mathsf{T}}} \nc{\sZ}{{\mathsf{Z}}} \nc{\sV}{{\mathsf{V}}}
\nc{\sfp}{{\mathsf{p}}} \nc{\sr}{{\mathsf{r}}}
\nc{\st}{{\mathsf{t}}} \nc{\sfb}{{\mathsf{b}}}
\nc{\sfc}{{\mathsf{c}}} \nc{\sd}{{\mathsf{d}}}
\nc{\sz}{{\mathsf{z}}}

\nc{\BK}{{\bar{K}}}

\nc{\tA}{{\widetilde{\mathbf{A}}}}
\nc{\tB}{{\widetilde{\mathcal{B}}}}
\nc{\tg}{{\widetilde{\mathfrak{g}}}} \nc{\tG}{{\widetilde{G}}}
\nc{\TM}{{\widetilde{\mathbb{M}}}{}}
\nc{\tO}{{\widetilde{\mathsf{O}}}{}}
\nc{\tU}{{\widetilde{\mathfrak{U}}}{}} \nc{\TZ}{{\tilde{Z}}}
\nc{\tx}{{\tilde{x}}} \nc{\tbv}{{\tilde{\bv}}}
\nc{\tfP}{{\widetilde{\mathfrak{P}}}{}} \nc{\tz}{{\tilde{\zeta}}}
\nc{\tmu}{{\tilde{\mu}}}

\nc{\urho}{\underline{\rho}} \nc{\uB}{\underline{B}}
\nc{\uC}{{\underline{\mathbb{C}}}} \nc{\ui}{\underline{i}}
\nc{\uj}{\underline{j}} \nc{\ofP}{{\overline{\mathfrak{P}}}}
\nc{\oB}{{\overline{\mathcal{B}}}}
\nc{\og}{{\overline{\mathfrak{g}}}} \nc{\oI}{{\overline{I}}}

\nc{\eps}{\varepsilon} \nc{\hrho}{{\hat{\rho}}}
\nc{\blambda}{{\boldsymbol{\lambda}}}

\nc{\one}{{\mathbf{1}}} \nc{\two}{{\mathbf{t}}}

\nc{\Rep}{{\mathop{\operatorname{\rm Rep}}}}
\nc{\Tot}{{\mathop{\operatorname{\rm Tot}}}}
\nc{\Ker}{{\mathop{\operatorname{\rm Ker}}}}
\nc{\Hilb}{{\mathop{\operatorname{\rm Hilb}}}}
\nc{\End}{{\mathop{\operatorname{\rm End}}}}
\nc{\Ext}{{\mathop{\operatorname{\rm Ext}}}}
\nc{\CHom}{{\mathop{\operatorname{{\mathcal{H}}\it om}}}}
\nc{\GL}{{\mathop{\operatorname{\rm GL}}}}
\nc{\gr}{{\mathop{\operatorname{\rm gr}}}}
\nc{\Id}{{\mathop{\operatorname{\rm Id}}}}
\nc{\defi}{{\mathop{\operatorname{\rm def}}}}
\nc{\length}{{\mathop{\operatorname{\rm length}}}}
\nc{\supp}{{\mathop{\operatorname{\rm supp}}}}

\nc{\Cliff}{{\mathsf{Cliff}}}
\nc{\Fl}{{\mathsf{Fl}}} \nc{\Fib}{{\mathsf{Fib}}}
\nc{\Coh}{{\mathsf{Coh}}} \nc{\FCoh}{{\mathsf{FCoh}}}

\nc{\reg}{{\text{\rm reg}}}

\nc{\cplus}{{\mathbf{C}_+}} \nc{\cminus}{{\mathbf{C}_-}}
\nc{\cthree}{{\mathbf{C}_*}} \nc{\Qbar}{{\bar{Q}}}

\nc{\bOmega}{{\overline{\Omega}}}

\nc{\seq}[1]{\stackrel{#1}{\sim}}

\nc{\aff}{\operatorname{aff}}

\begin{document}

\title{Equivariant $K$-theory of Hilbert
schemes via Shuffle Algebra}

\author{Boris Feigin}
 \address{Independent University of Moscow, 11 Bol'shoy Vlas'evskiy per., Moscow 119002,
 Russia AND \ \
 Landau Institute for Theoretical Physics, 1A pr. Akademika Semenova, Chernogolovka, Russia}
 \curraddr{Higher School of Economics, 20 Myasnitskaya ul., Moscow 101000, Russia}
 \email{borfeigin@gmail.com}

\author[Alexander Tsymbaliuk]{Alexander Tsymbaliuk}
 \address{Independent University of Moscow, 11 Bol'shoy Vlas'evskiy per., Moscow 119002, Russia}
 \curraddr{Department of Mathematics, MIT, 77 Massachusetts Ave., Cambridge, MA  02139, USA}
 \email{sasha\_ts@mit.edu}

\begin{abstract}
In this paper we construct actions of \textit{Ding-Iohara} and
\textit{shuffle} algebras on the sum of localized equivariant
$K$-groups of Hilbert schemes of points on $\mathbb C^2$. We show
that commutative elements $K_i$ of shuffle algebra act through
vertex operators over the positive part $\{\mathfrak{h}_i\}_{i>0}$
of the Heisenberg algebra in these $K$-groups. This provides an
action of the Heisenberg algebra itself. Finally, we normalize the
basis of the structure sheaves of fixed points in such a way that it
corresponds to the basis of Macdonald polynomials in the Fock space
$\mathbb C[\mathfrak{h}_1, \mathfrak{h}_2, \ldots]$.
\end{abstract}

\maketitle

\section{Introduction}
 For any surface $X$, let $X^{[n]}$ denote the Hilbert scheme of $n$
points on $X$. The Heisenberg algebra $\{ \mathfrak{h}_i \}_{i\in
\mathbb Z \backslash 0}$ is known (see~\cite{nak}) to act through
natural correspondences on the sum of cohomology rings
$\bigoplus_{n}\mathbb{H}^{*}(X^{[n]})$.

 From now on we deal only with the case $X=\mathbb C^2$.
Then one can consider localized equivariant cohomologies instead of
the usual ones. Let $R=\bigoplus_{n}\mathbb{H}^{2n}_{\mathbb T}
(X^{[n]}) \otimes_{\mathbb{H}_{\mathbb T}(\textrm{pt})}\on{Frac}
(\mathbb{H}_{\mathbb T}(\textrm{pt}))$.
 As shown in~\cite{l}, $R$ is isomorphic to the Fock space $\Lambda_F:=\mathbb
C(\hbar, \hbar')[\mathfrak{h}_1, \mathfrak{h}_2,\ldots]$, and
after certain normalization, there is an isomorphism $\Delta: R\to
\Lambda_F$ sending the basis of fixed points to Jack polynomials
and $\{\mathfrak{h}_i\}_{i>0}$ to operators of multiplication by
$p_i$.

 In this paper we construct an action of $1+\sum_{i>0}{\widetilde{K_i}z^i}:=\exp \left(\sum_{i>0}{((-1)^{i-1}/i)\mathfrak{h}_iz^i}\right)$
on the sum of localized equivariant $K$-groups
$M=\bigoplus_{n}K^{\mathbb T} (X^{[n]}) \otimes_{K^{\mathbb
T}(\textrm{pt})}\on{Frac} (K^{\mathbb T}(\textrm{pt}))$ in geometric
terms. This also provides the Heisenberg algebra's action. We find
an isomorphism $\Theta:M\to \Lambda_F$ which takes the normalized
fixed point basis $\{\langle\lambda\rangle\}$ to Macdonald
polynomials $\{P_{\lambda}\}$. Isomorphism $\Theta$ takes operators
$\widetilde{K_i}$ to operators of multiplication by $e_i$ acting in
the Macdonald polynomials basis via the Pieri formulas.

 For achieving this result and for its own sake, we construct representations of two algebras: $A$ and $S$ (called \textit{Ding-Iohara}
and \textit{shuffle} algebras correspondingly) on $M$. In fact,
subalgebra $\mathcal{S}$ of the \textit{shuffle} algebra $S$,
generated by $S_1$, is of particular interest to us. It is a
trigonometric analogue of the Feigin-Odesskii algebra, studied
in~\cite{fo}.

 The same operators appear in~\cite{sv}, where an action of the Hall
algebra of an elliptic curve is constructed on $R$. We strongly
recommend~\cite{sv} to an interested reader as its results are more
complete than ours and contain a completely different viewpoint. A
new insight on the relation between these two approaches was
recently obtained in~\cite{n}.

 In Section~\ref{section2}, we define Ding-Iohara
and shuffle algebras and remind some properties of them. In
Section~\ref{section3}, we construct an action of Ding-Iohara
algebra $A$ on $M$. In Section~\ref{section4}, we verify that these
operators do give a representation of $A$. In
Section~\ref{section5}, we define an action of shuffle algebra on
$M$. In Section~\ref{section6}, we present operators
$\widetilde{K_i}$, normalization of the fixed point basis $\{
\langle\lambda\rangle\}$, and an isomorphism $\Theta:M\iso
\Lambda_F$ with the aforementioned properties. Finally, in
Section~\ref{section7}, we consider a vector $v=\sum_{n\geq
0}{\left[\mathcal{O}_{X^{[n]}}\right]}$ (and its generalized version
$w(x)$).
We show that this vector is an eigenvector for the negative half of
our Heisenberg algebra. Finally, we determine the bilinear symmetric
form on $M$, for which operators $e_i$ and $f_i$ are adjoint to each
other.

\subsection*{Acknowledgments:} We are grateful to M.~Finkelberg, L.~Rybnikov, and
O.~Schiffmann for discussion and correspondence. We would like to
thank the referee for many useful suggestions.

\section{Ding-Iohara and shuffle algebras}
\label{section2}
  Let us fix three parameters $q_1,\ q_2,\ q_3$. Now we define the \textit{Ding-Iohara algebra} A. This is an associative algebra generated by $e_i,\ f_i,\
\psi^{\pm}_{j}\ (i \in \mathbb Z, j \in \mathbb Z_{+})$ with the
following defining relations:
\begin{equation}
\label{1}
e(z)e(w)(z-q_1w)(z-q_2w)(z-q_3w)=-e(w)e(z)(w-q_1z)(w-q_2z)(w-q_3z)
\end{equation}
\begin{equation}
\label{2}
f(z)f(w)(w-q_1z)(w-q_2z)(w-q_3z)=-f(w)f(z)(z-q_1w)(z-q_2w)(z-q_3w)
\end{equation}
\begin{equation}
\label{3} [e(z),
f(w)]=\frac{\delta(z/w)}{(1-q_1)(1-q_2)(1-q_3)}(\psi^{+}(w)-\psi^{-}(z))
\end{equation}
\begin{equation}
\label{4}
\psi^{\pm}(z)e(w)(z-q_1w)(z-q_2w)(z-q_3w)=-e(w)\psi^{\pm}(z)(w-q_1z)(w-q_2z)(w-q_3z)
\end{equation}
\begin{equation}
\label{5}
\psi^{\pm}(z)f(w)(w-q_1z)(w-q_2z)(w-q_3z)=-f(w)\psi^{\pm}(z)(z-q_1w)(z-q_2w)(z-q_3w)
\end{equation}
where these generating series are defined as follows:
$$e(z)=\sum_{i=-\infty}^{\infty}{e_iz^{-i}},\  \  f(z)=\sum_{i=-\infty}^{\infty}{f_iz^{-i}},\ \ \psi^{\pm}(z)=\sum_{j\geq 0}{\psi^{\pm}_{j}z^{\mp
j}},\ \ \delta(z)=\sum_{i=-\infty}^{\infty}{z^i}. $$

\begin{remark}

 \noindent{These relations are very similar to the relations of quantum
affine algebras (except for Serre relations).}
\end{remark}

 \noindent{We denote by $A_{+}$ (resp. $A_{-}$) the subalgebra of $A$ generated by
$e_i$ (resp. $f_i)$.}

\medskip \noindent
Next, following~\cite{fo}, we define the shuffle algebra $S$
depending on $q_1,\ q_2,\ q_3$. Fix a function
$$\lambda(x,y)=\frac{(x-q_1y)(x-q_2y)(x-q_3y)}{(x-y)^3}.$$
 Algebra $S$ is an associative graded algebra $S=\bigoplus_{n\geq
0}{S_n}$. Each graded component $S_n$ consists of rational functions
of the form $F(x_1,\ldots,x_n)=\frac{f(x_1,\ldots,x_n)}{\prod_{1\leq
i<j\leq n}{(x_i-x_j)^2}},$ where $f(x_1,\ldots,x_n)$ is a symmetric
Laurent polynomial. For $F\in S_m$ and $G\in S_n$, the product
$F*G\in S_{n+m}$ is defined by the formula
\begin{multline}\label{star-product}
(F*G)(x_1,\ldots,
x_{n+m})\\=\textrm{Sym}\left(F(x_1,\ldots,x_m)G(x_{m+1},\ldots,x_{m+n})\prod_{1\leq
i\leq m<j\leq m+n}{\lambda(x_i,x_j)}\right)
\end{multline}
 with $\textrm{Sym}$ standing for a symmetrization. It endows
$S$ with a structure of an associative algebra.

Now we formulate some known properties of shuffle algebras.

\begin{theorem}
\label{feigin1}
 (a) For any $q_1,q_2,q_3$, there is a natural homomorphism $\Xi:A_{+}\to S$, which takes $e_a \in A_{+}$ into $x^a\in S$.

 (b) For generic parameters $q_1,\ q_2,\ q_3$, homomorphism $\Xi$ is an isomorphism; in particular, the whole algebra $S$ is generated by $S_1$.
\end{theorem}

\begin{theorem}
\label{feigin3} \footnote{\ This theorem was conjectured by the
first author and its proof follows from the recent
papers~\cite{n},~\cite{sch}.}
 For generic $q_1,\ q_2$ and $q_3:=q_1^{-1}q_2^{-1}$, the subalgebra $\mathcal{S}$
generated by $S_1$ consists of all rational functions of the form
$F(x_1,\ldots,x_n)=f(x_1,\ldots,x_n)\prod_{1\leq i<j\leq
n}{(x_i-x_j)^{-2}},$ where $f(x_1,\ldots,x_n)$ is a symmetric
Laurent polynomial satisfying 'wheel conditions', that is
$f(x_1,\ldots,x_n)=0$ if $x_1/x_2=q_1$ and $x_2/x_3=q_j$ for $j=2,
3$.
\end{theorem}

\emph{Remark:}
 For arbitrary $q_1,q_2,q_3$ the connection between
the Ding-Iohara algebra $A$ and the shuffle subalgebra $\mathcal{S}$
is as follows. Let $I_e$ be the kernel of $\Xi$ from
Theorem~\ref{feigin1}(a) and $I_f$ be the transposed ideal of
$A_{-}$, that is, $I_f$ is obtained from $I_e$ via the involution
$e_i\longmapsto f_{-i}$. Then the factor of $A$ by the ideals $I_f,
I_e$ is what we are most interested in. It may be viewed as a
Drinfeld double of the shuffle subalgebra $\mathcal{S}$
(see~\cite{n},~\cite{sv} for more details). In the setup of
Theorem~\ref{feigin3}, the description of the ideal $I_e$ had been
conjectured in~\cite{ffjmm} and was recently proved in~\cite{sch}.

\begin{theorem}
\label{feigin2}

For each $n\geq 1$, define elements $K_n\in S_n$ by
$$K_2(x_1,x_2):=\frac{(x_1-q_1x_2)(x_2-q_1x_1)}{(x_1-x_2)^2},\ \ \
\ \ K_n(x_1,\ldots,x_n):=\prod_{1\leq i<j\leq n}{K_2(x_i,x_j)}.$$ In
particular, $K_1(x_1)=1$. If $q_1q_2q_3=1$, elements $K_n\in S_n$
commute.
\end{theorem}

 \noindent{The subalgebra generated by $K_i$ is studied in~\cite{fhsy}; in particular, Theorem~\ref{feigin2} is proved therein.}

\medskip

 Our work was motivated by~\cite{ffnr} and~\cite{ts}.

\section{Construction of operators}
\label{section3}

\subsection{Correspondences}
We recall that $X=\mathbb C^2$ throughout this paper. In this case
the Hilbert scheme of $n$ points $X^{[n]}$ is identified
(set-theoretically) with the set of all codimension $n$ ideals of
$\BC[x,y]$. Let us recall correspondences used by H.~Nakajima to
construct a representation of the Heisenberg algebra on
$\bigoplus_{n}\mathbb H^{*}(X^{[n]})$. This action is constructed
via correspondences $P[i]\subset \coprod_{n}{X^{[n]}\times
X^{[n+i]}}$. Let us mention their definition for an arbitrary $i$,
though we will actually need only $P[1]$ and $P[-1]$. For $i>0$, the
correspondence $P[i]\subset \coprod_{n}{X^{[n]}\times X^{[n+i]}}$
consists of all pairs of ideals $(J_1, J_2)$ of $\BC[x,y]$ of
codimension $n,\ n+i$ respectively, such that $J_2\subset J_1$ and
the factor $J_1/J_2$ is supported at a single point (the latter
condition is automatic for $i=1$). For $i<0$, $P[i]$ is defined as
the transposed to $P[-i]$. Let $L$ be a tautological line bundle on
$P[1]$ whose fiber at any point $(J_1, J_2) \in P[1]$ equals
$J_1/J_2$. There are natural projections $\bp, \bq$ from $P[1]$ to
$X^{[n]}$ and $X^{[n+1]}$, correspondingly.

\subsection{Fixed points}
\label{fixed points} There is a natural action of $\mathbb T=\mathbb
C^{*}\times \mathbb C^{*}$ on each $X^{[n]}$ induced from the one on
$X$ given by the formula $(t_1,t_2)(x,y)=(t_1\cdot x, t_2\cdot y)$.
Set $(X^{[n]})^{\mathbb T}$ of $\mathbb T$-fixed points in $X^{[n]}$
is finite, and is in bijection with size $n$ Young diagrams. Namely,
for each Young diagram $\lambda=(\lambda_1,\ldots, \lambda_k)$ we
have an ideal
 $\BC[x,y]\cdot (\BC x^{\lambda_1}y^0\oplus\cdots\oplus\BC x^{\lambda_k}y^{k-1}\oplus\BC y^k)=:J_{\lambda}\in (X^{[n]})^{\mathbb T}$.

\subsection{Equivariant $K$-groups}
We denote the direct sum of equivariant (complexified) $K-$groups as
${}'M:=\bigoplus_{n}K^{\mathbb T}(X^{[n]}).$ It is a module over
$K^{\mathbb T}(\textrm{pt})=\BC[\mathbb T]= \BC[t_1,t_2]$. We define
$$M:=\ {}'M\otimes_{K^{\mathbb T}(\textrm{pt})} \on{Frac}(K^{\mathbb
T}(\textrm{pt}))={}'M\otimes_{\BC [t_1,t_2]} \BC (t_1,t_2).$$

It has a natural grading: $M=\bigoplus_{n}M_{n},\ M_{n}=K^{\mathbb
T}(X^{[n]}) \otimes_{K^{\mathbb T}(\textrm{pt})}
\on{Frac}(K^{\mathbb T}(\textrm{pt})).$

According to the Thomason localization theorem, restriction to the
$\mathbb T$-fixed point set induces an isomorphism
$$K^{\mathbb T}(X^{[n]})
\otimes_{K^{\mathbb T}(\textrm{pt})} \on{Frac}(K^{\mathbb
T}(\textrm{pt}))\iso K^{\mathbb T}((X^{[n]})^{\mathbb T})
\otimes_{K^{\mathbb T}(\textrm{pt})} \on{Frac}(K^{\mathbb
T}(\textrm{pt})).$$

The structure sheaves $\{\lambda\}$ of the $\mathbb T$-fixed points
$J_\lambda$ (see Section~\ref{fixed points}) form a basis in
$\bigoplus_{n}K^{\mathbb T} ((X^{[n]})^{\mathbb T})
\otimes_{K^{\mathbb T}(\textrm{pt})}\on{Frac} (K^{\mathbb
T}(\textrm{pt}))$. Since embedding of a point $J_\lambda$ into
$X^{[n]}$ is a proper morphism, the direct image in the equivariant
$K$-theory is well defined, and we denote by $[\lambda]\in M_n$ the
direct image of the structure sheaf $\{\lambda\}$. The set
$\{[\lambda]\}$ forms a basis of $M$.

\subsection{Representation of Ding-Iohara algebra on M}

Let us now consider the tautological vector bundle $\mathfrak{F}$ on
$X^{[n]}$, whose fiber at the point corresponding to an ideal $J$
equals $\BC[x,y]/J$. We introduce generating series $\ba(z),\
\bc(z)$ (with coefficients in $M$) as follows:
$$\ba(z):=\Lambda^{\bullet}_{-1/z}(\mathfrak{F})=\sum_{i\geq
0}{[\Lambda^i(\mathfrak{F})](-1/z)^i},$$
$$\bc(z):=\ba(zt_1)\ba(zt_2)\ba(zt_1^{-1}t_2^{-1})\ba(zt_1^{-1})^{-1}\ba(zt_2^{-1})^{-1}\ba(zt_1t_2)^{-1}.$$

We also define the operators
\begin{equation}
\label{raz} e_i=\bq_*(L^{\otimes i}\otimes \bp^*):\ M_{n}\to M_{n+1}
\end{equation}
\begin{equation}
\label{dva} f_i=\bp_*(L^{\otimes (i-1)}\otimes \bq^*):\ M_{n}\to
M_{n-1}
\end{equation}

 So $e_i$ is a composition of pulling back along $P[1]\to M_n$, tensoring
by $L^{\otimes i}$, and finally pushing forward along $P[1]\to
M_{n+1}$, while $f_{i+1}$ is obtained by the inverse order of
these operations.

 We consider the following generating series of operators acting on $M$:
\begin{equation}
\label{dvas}
 e(z)=\sum_{r=-\infty}^\infty e_rz^{-r}:\ M_{n}\to M_{n+1}[[z,z^{-1}]]
\end{equation}
\begin{equation}
\label{tris}
 f(z)=\sum_{r=-\infty}^\infty f_rz^{-r}:\ M_{n}\to M_{n-1}[[z,z^{-1}]]
\end{equation}
\begin{equation}
\label{raz1}
 \psi^{\pm}(z)\mid_{M_{n}}=\sum_{r=0}^{\infty} \psi^{\pm}_rz^{\mp r}:=\left(-\frac{1-t_1^{-1}t_2^{-1}z^{-1}}{1-z^{-1}}\bc(z)\right)^{\pm} \in M_{n}[[z^{\mp 1}]]
\end{equation}
 where $\left(\ldots\right)^{\pm}$ denotes the expansion in $z^{\mp 1}$, respectively.

Formula~(\ref{raz1}) should be understood as follows:
$\psi^{\pm}(z)$ acts by multiplication in $K$-theory by
$\left(-\frac{1-t_1^{-1}t_2^{-1}z^{-1}}{1-z^{-1}}\bc(z)\right)^{\pm}$,
and $\psi_r^{\pm}$ are defined as the coefficients of these series.

\begin{theorem}
\label{main}

 Operators $e_i, f_i,\psi^{\pm}_j$, defined in (\ref{dvas}--\ref{raz1}), satisfy
 relations~(\ref{1}--\ref{5}) with parameters $q_1=t_1,\ q_2=t_2,\
 q_3=t_1^{-1}t_2^{-1}$; that is, they give rise to a representation of the algebra $A$ on $M$.
\end{theorem}

 We prove this theorem in Section 4, by a straightforward
 verification of all relations.
\medskip

\emph{Notation:}  For any Young diagram $\lambda=(\lambda_1\geq
\cdots\geq\lambda_k)$ and its box $\square_{i,j}$ with the
coordinates $(i,j)$ (i.e., it stands in the $i$th row and $j$th
column), where $1\leq i\leq k,\ 1\leq j\leq \lambda_i$, we introduce
functions $l(\square)$ and $a(\square)$, called \textit{legs} and
\textit{arms}, correspondingly:
$$l(\square):=\lambda_i-j,\ \ \ \ \ \  a(\square):=\max\{k|\lambda_k\geq j\}-i.$$
We also denote by $\Sigma_{1}(\square)$ all boxes of $\lambda$ with
the coordinates $(i,k<j)$ and by $\Sigma_{2}(\square)$ all boxes of
$\lambda$ with the coordinates $(k<i,j)$. Sometimes we write
$\lambda+j$ for the diagram $\lambda+\square_{j,\lambda_j+1}$ if it
makes sense (i.e., if it is still a diagram). Finally, we call box
$\square_{i,j}$ a \textit{corner} if $j=\lambda_i>\lambda_{i+1}$ and
a \textit{hole} if $j=\lambda_i+1\leq \lambda_{i-1}.$

Now we compute the matrix coefficients of operators $e_i, f_i$ and
the eigenvalues of $\psi^{\pm}(z)$ in the fixed point basis.

\begin{lemma}
\label{eigenv}

\textrm{(a)} The only nonzero matrix coefficients of the operators
$e_i, f_i$ in the fixed point basis $[\lambda]$ of $M$ are as
follows:
\begin{equation*}
 e_{i[\lambda,\lambda+k]}=(1-t_1)^{-1}(1-t_2)^{-1}(t_1^{\lambda_k}t_2^{k-1})^i
\prod_{s\in
\Sigma_1(\square_{k,\lambda_k+1})}{\frac{1-t_1^{-l(s)+1}t_2^{a(s)+1}}{1-t_1^{-l(s)}t_2^{a(s)+1}}}
\prod_{s\in
\Sigma_2(\square_{k,\lambda_k+1})}{\frac{1-t_1^{l(s)+1}t_2^{-a(s)+1}}{1-t_1^{l(s)+1}t_2^{-a(s)}}}
\end{equation*}
\begin{equation*}
f_{i[\lambda,\lambda-k]}=(t_1^{\lambda_k-1}t_2^{k-1})^{i-1}\prod_{s\in
\Sigma_1(\square_{k,\lambda_k})}{\frac{1-t_1^{l(s)+1}t_2^{-a(s)}}{1-t_1^{l(s)}t_2^{-a(s)}}}
\prod_{s\in
\Sigma_2(\square_{k,\lambda_k})}{\frac{1-t_1^{-l(s)}t_2^{a(s)+1}}{1-t_1^{-l(s)}t_2^{a(s)}}}
\end{equation*}

\textrm{(b)} The eigenvalue of $\psi^{\pm}(z)$ on $[\lambda]$
equals
$$\left(-\frac{1-t_1^{-1}t_2^{-1}z^{-1}}{1-z^{-1}}
\prod_{\square \in
\lambda}{\frac{(1-t_1^{-1}\chi(\square)z^{-1})(1-t_2^{-1}\chi(\square)z^{-1})(1-t_1t_2\chi(\square)z^{-1})}
{(1-t_1\chi(\square)z^{-1})(1-t_2\chi(\square)z^{-1})(1-t_1^{-1}t_2^{-1}\chi(\square)z^{-1})}}\right)^{\pm}$$
where $\chi(\square_{i,j})=t_1^{j-1}t_2^{i-1}.$
\end{lemma}

\begin{proof}

(a) For $(\lambda, \lambda')\in P[1]$, let
$\rho:J_{\lambda'}\hookrightarrow J_{\lambda},\
\pi:k[x,y]/J_{\lambda'} \twoheadrightarrow k[x,y]/J_{\lambda}$ be
the natural maps. The tangent space
$\mathfrak{T}_{(J_{\lambda},J_{\lambda'})}(P[1])$ is a kernel of
the map
$$\textrm{Hom}(J_{\lambda'},k[x,y]/J_{\lambda'})\oplus
\textrm{Hom}(J_{\lambda},k[x,y]/J_{\lambda})\twoheadrightarrow
\textrm{Hom}(J_{\lambda'},k[x,y]/J_{\lambda}),\ \ \ (\alpha, \beta)
\longmapsto \pi \circ \alpha - \beta \circ \rho.$$
 Further, we write simply $\lambda$ instead of $J_{\lambda}$.

Let us denote by $\chi_{(\lambda,\lambda')}$ the character of
$\mathbb T$ in the tangent space
$\mathfrak{T}_{(\lambda,\lambda{}')}(P[1])$ and by
${\chi(L)}_{(\lambda,\lambda')}$ the character of $\mathbb T$ in
the fiber of $L$ at the point $(\lambda,\lambda{}')$. We write
$S\chi_{(\lambda)}$ (resp., $S\chi_{(\lambda,\lambda{}')}$) for
the character of $\mathbb T$ in the symmetric algebra
$\textrm{Sym}^{\bullet}\mathfrak T_{(\lambda)}X^{[n]}$ (resp.,
$\textrm{Sym}^{\bullet}\mathfrak T_{(\lambda,\lambda{}')}P[1]$).

According to the Bott-Lefschetz fixed point formula, the matrix
coefficient $\bp_*(L^{\otimes i}\otimes \bq^*)_{[\lambda',\lambda]}$
of $\bp_*(L^{\otimes i}\otimes \bq^*):\ M_{n+1}\to M_{n}$ with
respect to the basis elements $[\lambda]\in K^{\mathbb T}
(X^{[n]}),\ [\lambda']\in K^{\mathbb T} (X^{[n+1]})$ equals
$\chi(L)^i_{(\lambda,\lambda')}S\chi_{(\lambda,\lambda')}/
S\chi_{(\lambda')}$. Similarly, the matrix coefficient
$\bq_*(L^{\otimes i}\otimes \bp^*)_{(\lambda,\lambda')}$ of
$\bq_*(L^{\otimes i}\otimes \bp^*):\ M_{n}\to M_{n+1}$ with respect
to the basis elements $[\lambda]\in K^{\mathbb T} (X^{[n]}),\
[\lambda']\in K^{\mathbb T} (X^{[n+1]})$ equals
$\chi(L)^i_{(\lambda,\lambda')}S\chi_{(\lambda,\lambda')}/
S\chi_{(\lambda)}$.

 Now it is straightforward to check the formulas.

(b) This follows from the multiplicativity of
$\Lambda^{\bullet}_z(F)$ on the long exact sequences and the fact
that $\{\chi(\square)|\square \in \lambda\}$ is a set of $\mathbb
T$--characters at the fiber $\mathfrak{F}|_{\lambda}$.
\end{proof}

 It will be convenient for us to have more uniform formulas for the matrix coefficients of $e_i$ and $f_i$.
 Those are provided by the following proposition:

\begin{proposition}
\label{eig'} For any Young diagram $\lambda$ and $i\in \BN$ we have
$$e_{r[\lambda-i,\lambda]}=\frac{(t_1^{\lambda_i-1}t_2^{i-1})^r}{(1-t_1^{\lambda_1-\lambda_i+1}t_2^{1-i})(1-t_1t_2)}
\prod_{j=1}^{\infty}{\frac{1-t_1^{\lambda_j-\lambda_i+1}t_2^{j-i+1}}{1-t_1^{\lambda_{j+1}-\lambda_i+1}t_2^{j-i+1}}},$$
$$f_{r[\lambda+i,\lambda]}=\frac{(t_1^{\lambda_i}t_2^{i-1})^{r-1}(1-t_1^{\lambda_i-\lambda_1+1}t_2^{i})}{1-t_1t_2}
\prod_{j=1}^{\infty}{\frac{1-t_1^{\lambda_i-\lambda_{j+1}+1}t_2^{i-j}}{1-t_1^{\lambda_{i}-\lambda_j+1}t_2^{i-j}}}.$$
\end{proposition}

\begin{proof}
 These are formulas of Lemma~\ref{eigenv}(a), written in a different way.
\end{proof}

\medskip

\section{Proof of Theorem~\ref{main}}
\label{section4}

\begin{definition}

\textsl{We denote the elementary symmetric polynomials in $q_1, q_2,
q_3$ by $\sigma_1,  \sigma_2, \sigma_3$:
$$\sigma_1:=q_1+q_2+q_3=t_1+t_2+t^{-1}_1t^{-1}_2,\
\sigma_2:=q_1q_2+q_1q_3+q_2q_3=t^{-1}_1+t^{-1}_2+t_1t_2,\
\sigma_3:=q_1q_2q_3=1.$$}
\end{definition}

\textbf{\emph{Convention:}} In this section we check
(\ref{1}--\ref{5}) explicitly in the fixed point basis. While
comparing expressions of left-hand side and right-hand side, we
denote by $P_i$ the mutual factor.

\medskip

First, let us check~(\ref{1}).

\begin{proof}

  For any integers $i, j$ we have to prove the following equation:
$$e_{i+3}e_j-\sigma_1e_{i+2}e_{j+1}+\sigma_2e_{i+1}e_{j+2}-\sigma_3e_ie_{j+3}$$$$=\sigma_3e_je_{i+3}-\sigma_2e_{j+1}e_{i+2}+\sigma_1e_{j+2}e_{i+1}-e_{j+3}e_i.$$
Let us compare the matrix coefficients of the left-hand side and the
right-hand side on any pair of Young diagrams $[\lambda,
\lambda'=\lambda+\square_{i_1,j_1}+\square_{i_2,j_2}]$.

\noindent(a) Suppose $i_1=i_2$; that is, the added two boxes lie in
the same row. Then:
$$(e_{i+3}e_j-\sigma_1e_{i+2}e_{j+1}+\sigma_2e_{i+1}e_{j+2}-\sigma_3e_ie_{j+3})_{[\lambda,\lambda{}']}=
(\ldots)(1-\sigma_1t_1^{-1}+\sigma_2t_1^{-2}-\sigma_3t_1^{-3})=0$$
since $t_1^{-1}$ is a root of polynomial
$1-\sigma_1t+\sigma_2t^2-\sigma_3t^3$.

Similarly:
$(\sigma_3e_je_{i+3}-\sigma_2e_{j+1}e_{i+2}+\sigma_1e_{j+2}e_{i+1}-e_{j+3}e_i)_{[\lambda,\lambda{}']}=0$.

\noindent(b) Suppose $j_1=j_2$; that is, the added two boxes lie in
the same column.

 This case is entirely similar since $t_2^{-1}$ is also a root of $1-\sigma_1t+\sigma_2t^2-\sigma_3t^3.$

\noindent(c) Suppose $i_1<i_2, j_1>j_2$.

 According to the formulas of Lemma~\ref{eigenv}(a), the only difference occurs at the box $\square_{i_1,j_2}$.

 Let us denote $a:=j_1-j_2,\ b:=i_2-i_1,\ \chi_1:=t_1^{j_1-1}t_2^{i_1-1},\ \chi_2:=t_1^{j_2-1}t_2^{i_2-1}$.
Then
$$(e_{i+3}e_j-\sigma_1e_{i+2}e_{j+1}+\sigma_2e_{i+1}e_{j+2}-\sigma_3e_ie_{j+3})_{[\lambda,\lambda{}']}=$$$$
P_1(1-t_1^{-a}t_2^b)^{-1}(1-t_1^{-a+1}t_2^b)(1-t_1^{a+1}t_2^{-b})^{-1}(1-t_1^{a+1}t_2^{-b+1})
\chi_1^{j}\chi_2^{i+3}
\left(1-\sigma_1\left(\frac{\chi_1}{\chi_2}\right)+\sigma_2\left(\frac{\chi_1}{\chi_2}\right)^2-\sigma_3\left(\frac{\chi_1}{\chi_2}\right)^3\right)$$
$$+P_1(1-t_1^{a}t_2^{-b})^{-1}(1-t_1^{a}t_2^{-b+1})(1-t_1^{-a}t_2^{b+1})^{-1}(1-t_1^{-a+1}t_2^{b+1})
\chi_1^{i+3}\chi_2^{j}\left(1-\sigma_1\left(\frac{\chi_2}{\chi_1}\right)+\sigma_2\left(\frac{\chi_2}{\chi_1}\right)^2-\sigma_3\left(\frac{\chi_2}{\chi_1}\right)^3\right),
$$

$$(\sigma_3e_je_{i+3}-\sigma_2e_{j+1}e_{i+2}+\sigma_1e_{j+2}e_{i+1}-e_{j+3}e_i)_{[\lambda,\lambda{}']}=$$$$
P_1(1-t_1^{a}t_2^{-b})^{-1}(1-t_1^{a}t_2^{-b+1})(1-t_1^{-a}t_2^{b+1})^{-1}(1-t_1^{-a+1}t_2^{b+1})
\chi_1^{j}\chi_2^{i+3}\left(\sigma_3-\sigma_2\left(\frac{\chi_1}{\chi_2}\right)+\sigma_1\left(\frac{\chi_1}{\chi_2}\right)^2-\left(\frac{\chi_1}{\chi_2}\right)^3\right)$$
$$+P_1(1-t_1^{-a}t_2^b)^{-1}(1-t_1^{-a+1}t_2^b)(1-t_1^{a+1}t_2^{-b})^{-1}(1-t_1^{a+1}t_2^{-b+1})
\chi_1^{i+3}\chi_2^{j}\left(\sigma_3-\sigma_2\left(\frac{\chi_2}{\chi_1}\right)+\sigma_1\left(\frac{\chi_2}{\chi_1}\right)^2-\left(\frac{\chi_2}{\chi_1}\right)^3\right).
$$

 Set $u:=t_1^a,\ v:=t_2^b$, so that $\chi_2/\chi_1=u^{-1}v$. The first
summand of the left-hand side equals
$$P_1\chi_1^{j}\chi_2^{i+3}(u-v)^{-1}(u-t_1v)(v-t_1u)^{-1}(v-t_1t_2u)(v-t_1u)(v-t_2u)(v-t_1^{-1}t_2^{-1}u)v^{-3}$$$$=
P_1\chi_1^{j}\chi_2^{i+3}(u-v)^{-1}(u-t_1v)(v-t_1t_2u)(v-t_2u)(v-t_1^{-1}t_2^{-1}u)v^{-3},$$
while the first summand of the right-hand side equals
$$P_1\chi_1^j\chi_2^{i+3}(v-u)^{-1}(v-t_2u)(u-t_2v)^{-1}(u-t_1t_2v)(u-t_1v)(u-t_2v)(u-t_1^{-1}t_2^{-1}v)(-u^3v^{-3})u^{-3}$$
$$=P_1\chi_1^j\chi_2^{i+3}(u-v)^{-1}(v-t_2u)(u-t_1t_2v)(u-t_1v)(u-t_1^{-1}t_2^{-1}v)v^{-3}.$$
 These two expressions coincide. In the same way we check the
equality of the second summands. This completes the proof in this
case.

\noindent(d) Suppose $i_1>i_2, j_1<j_2$. This case follows from (c).
\end{proof}

\noindent  Equation~(\ref{2}) is entirely similar to the one above,
so we omit
 it.

\medskip

  Now we compute $[e(z), f(w)]$. We prove the following
proposition at first.
\begin{proposition}
\label{diagonal}

 Coefficients of the the series $[e(z), f(w)]$ are diagonalizable
in the fixed point basis $[\lambda]$.
\end{proposition}

\begin{proof}

 We have to verify $(e_if_j)_{[\lambda,\lambda{}']}=(f_je_i)_{[\lambda,\lambda{}']}$ for any
 diagrams  $[\lambda,\lambda'=\lambda+\square_{i_1,j_1}-\square_{i_2,j_2}]$ with $(i_1,j_1)\ne (i_2,j_2)$.

Let us consider the case $i_1<i_2,\ j_1>j_2$ (the case $i_1>i_2,\
j_1<j_2$ is completely analogous). We define $a:=j_1-j_2,\
b:=i_2-i_1.$ Then

$$(e_if_j)_{[\lambda,\lambda{}']}=P_2(1-t_1^{1-a}t_2^b)^{-1}(1-t_1^{1-a}t_2^{b+1})(1-t_1^{-a}t_2^b)^{-1}(1-t_1^{1-a}t_2^b)$$$$=P_2(1-t_1^{1-a}t_2^{b+1})(1-t_1^{-a}t_2^b)^{-1},$$
$$(f_je_i)_{[\lambda,\lambda{}']}=P_2(1-t_1^{-a}t_2^{b+1})^{-1}(1-t_1^{1-a}t_2^{b+1})(1-t_1^{-a}t_2^b)^{-1}(1-t_1^{-a}t_2^{b+1})$$$$=P_2(1-t_1^{1-a}t_2^{b+1})(1-t_1^{-a}t_2^b)^{-1}.
$$

So
$(e_if_j)_{[\lambda,\lambda{}']}=(f_je_i)_{[\lambda,\lambda{}']}$.
\end{proof}

\medskip
\noindent
 Next, we introduce series of operators
$\Phi^{\pm}(z)=\sum_{i=0}^{\infty}{\phi^{\pm}_iz^{\mp i}}$,
diagonalizable in the fixed point basis and satisfying the equation
$$[e(z),f(w)]=\frac{\delta(z/w)}{(1-t_1)(1-t_2)(1-t_1^{-1}t_2^{-1})}(\Phi^{+}(w)-\Phi^{-}(z)).$$
 Actually, this determines $\phi^{\pm}_{> 0}$ uniquely. We also set $\phi^{+}_0=-1, \phi^{-}_0=-1/t_1t_2$. Next, we check
\begin{equation}
\label{4'}
\phi^{\pm}(z)e(w)(z-q_1w)(z-q_2w)(z-q_3w)=-e(w)\phi^{\pm}(z)(w-q_1z)(w-q_2z)(w-q_3z)
\end{equation}
\begin{equation}
\label{5'}
\phi^{\pm}(z)f(w)(w-q_1z)(w-q_2z)(w-q_3z)=-f(w)\phi^{\pm}(z)(z-q_1w)(z-q_2w)(z-q_3w)
\end{equation}

\noindent Finally, we prove $\psi^{\pm}_i=\phi^{\pm}_i$, and
so~(\ref{4'}--\ref{5'}) imply~(\ref{4}--\ref{5}). This will prove
Theorem~\ref{main}.

\medskip

 Firstly, Proposition~\ref{diagonal} and Lemma~\ref{eigenv}(a) imply that $[e(z),f(w)]$ is
diagonalizable in the fixed point basis, and its eigenvalue at
$[\lambda]$ equals
$$\sum_{a,b\in \mathbb Z}{z^{-a}w^{-b}\gamma_{a+b}},$$ where
$$\gamma_i=(1-t_1)^{-1}(1-t_2)^{-1} \sum_{\square-\textrm{corner}}\left(\prod_{s \in
\Sigma_1(\square)}{\left[\frac{(1-t_1^{l(s)+1}t_2^{-a(s)})(1-t_1^{-l(s)+1}t_2^{a(s)+1})}{(1-t_1^{l(s)}t_2^{-a(s)})(1-t_1^{-l(s)}t_2^{a(s)+1})}\right]}\right.$$
$$\times\left.\prod_{s\in
\Sigma_2(\square)}{\left[\frac{(1-t_1^{-l(s)}t_2^{a(s)+1})(1-t_1^{l(s)+1}t_2^{-a(s)+1})}{(1-t_1^{-l(s)}t_2^{a(s)})(1-t_1^{l(s)+1}t_2^{-a(s)})}
\right]}\right)\chi^{i-1}(\square)$$
$$-(1-t_1)^{-1}(1-t_2)^{-1}
\sum_{\square-\textrm{hole}}\left(\prod_{s \in
\Sigma_1(\square)}{\left[\frac{(1-t_1^{l(s)+1}t_2^{-a(s)})(1-t_1^{-l(s)+1}t_2^{a(s)+1})}{(1-t_1^{l(s)}t_2^{-a(s)})(1-t_1^{-l(s)}t_2^{a(s)+1})}\right]}\right.$$
$$\times\left.\prod_{s\in
\Sigma_2(\square)}{\left[\frac{(1-t_1^{-l(s)}t_2^{a(s)+1})(1-t_1^{l(s)+1}t_2^{-a(s)+1})}{(1-t_1^{-l(s)}t_2^{a(s)})(1-t_1^{l(s)+1}t_2^{-a(s)})}
\right]}\right)\chi^{i-1}(\square).$$

\medskip
So as we want an equality
$$[e(z),f(w)]=\frac{\delta\left(\frac{z}{w}\right)\left(\Phi^{+}(w)-\Phi^{-}(z)\right)}{(1-t_1)(1-t_2)(1-t_1^{-1}t_2^{-1})}=$$$$
\frac{\sum_{a+b>0}{z^{-a}w^{-b}\phi_{a+b}^{+}}-\sum_{a+b<0}{z^{-a}w^{-b}\phi_{-a-b}^{-}}+\sum_{a+b=0}{z^{-a}w^{-b}(\phi_{0}^{+}-\phi_{0}^{-})}}{(1-t_1)(1-t_2)(1-t_1^{-1}t_2^{-1})}$$
to hold, we determine $\phi_{>0}^{+},\ -\phi_{>0}^{-},\
\phi_{0}^{+}-\phi_{0}^{-} $ uniquely. Let us, finally specialize
values $\phi^{\pm}_0$.

\noindent
 Next lemma is crucial for this:

\begin{lemma}
\label{Heis1}

 We have $[e_0, f_0]\mid_{\lambda}=-\frac{1}{(1-t_1)(1-t_2)},\ [e_0, f_1]\mid_{\lambda}=-\frac{1}{(1-t_1)(1-t_2)}+\sum_{\square\in \lambda}{\chi(\square)}.$
\end{lemma}

\begin{corollary}

 The operator $[e_0,f_1-f_0]$ is the operator of multiplication by $\det(\mathfrak{F})$.
\end{corollary}

\begin{proof} [Proof of Lemma~\ref{Heis1}]

 We will need another expression for
$\gamma_s=[e_0,f_s]$, obtained by applying Proposition~\ref{eig'}
instead of Lemma~\ref{eigenv}(a). For any Young diagram
$\lambda=(\lambda_1,\lambda_2,\ldots)$, we define
$\chi_i:=t_1^{\lambda_i-1}t_2^{i-1}$. Note that
$\chi_i=t_1^{-1}t_2^{i-1}$ for all $i\gg1$ (the \emph{stabilizing
condition}). Then:
\begin{multline}
\label{KUKU} \gamma_s=(1-t_1)^{-2}
\sum_{i=1}^{l(\lambda)+1}{\chi_i^{s-1}\prod_{i\ne
j}^{\infty}{\frac{(1-t_1^{\lambda_i-\lambda_j}t_2^{i-j+1})(1-t_1^{\lambda_j-\lambda_i+1}t_2^{j-i+1})}
{(1-t_1^{\lambda_i-\lambda_j}t_2^{i-j})(1-t_1^{\lambda_j-\lambda_i+1}t_2^{j-i})}}}\\-
(1-t_1)^{-2}\sum_{i=1}^{l(\lambda)+1}{(t_1\chi_i)^{s-1}\prod_{i\ne
j}^{\infty}{\frac{(1-t_1^{\lambda_j-\lambda_i}t_2^{j-i+1})(1-t_1^{\lambda_i-\lambda_j+1}t_2^{i-j+1})}
{(1-t_1^{\lambda_j-\lambda_i}t_2^{j-i})(1-t_1^{\lambda_i-\lambda_j+1}t_2^{i-j})}}}
\end{multline}

 \noindent (a) First, we prove $[e_0, f_0]\mid_{\lambda}=-(1-t_1)^{-1}(1-t_2)^{-1}$ for any $\lambda$.

 This is obvious for an empty diagram (and straightforward for any
$1$ row diagrams $\lambda$). Thus, it suffices to prove that
$[e_0,f_0]\mid_{\lambda}$ is independent of $\lambda$.
 Choose $k\in \BN$, such that $\lambda_{k-1}=0$ (we do not claim $\lambda_{k-2}\ne 0$).
Then $\chi_i=t_1^{-1}t_2^{i-1}$ for $i\geq k-1$, and so~(\ref{KUKU})
implies:
$$\gamma_0=(1-t_1)^{-2}\left(
\sum_{i=1}^{k}{\chi_i^{-1}\frac{\chi_i(1-t_1t_2^{2-k}\chi_i)}{\chi_i-t_2^{k-1}}\prod_{1\leq
j\ne i}^{k-1}{\frac{(\chi_j-t_2\chi_i)(\chi_i-t_1t_2\chi_j)}
{(\chi_j-\chi_i)(\chi_i-t_1\chi_j)}}}\right.$$$$-\left.
\sum_{i=1}^{k}{(t_1\chi_i)^{-1}\frac{\chi_i(1-t_1^2t_2^{2-k}\chi_i)}{\chi_i-t_1^{-1}t_2^{k-1}}\prod_{1\leq
j\ne i}^{k-1}{\frac{(\chi_i-t_2\chi_j)(\chi_j-t_1t_2\chi_i)}
{(\chi_i-\chi_j)(\chi_j-t_1\chi_i)}}}\right).$$

 The right-hand side is a rational expression in $\chi_i\ (1\leq i\leq k-2)$. Moreover, the degree
of its numerator is not greater than that of the denominator, while
it might have poles only at $\chi_i=\chi_j,\ \chi_i=t_1\chi_j$ or at
$\chi_i=t_2^{k-1}, t_1^{-1}t_2^{k-1}$. However, in the latter two
cases, there are no poles (since those coming from the denominator
are compensated by zeros of $\chi_i-t_1t_2\chi_{k-1}$ or
$\chi_i-t_2\chi_{k-1}$ coming from the numerator). The first two
poles (at $\chi_i=\chi_j,\ \chi_i=t_1\chi_j$) are simple, but a
straightforward computation verifies that residues at these points
are in fact zero. Hence, this rational function is constant. This
completes the proof of $[e_0,f_0]=-(1-t_1)^{-1}(1-t_2)^{-1}$.

\medskip
\noindent (b) Proof of $[e_0,
f_1]\mid_{\lambda}=-(1-t_1)^{-1}(1-t_2)^{-1}+\sum_{\square\in
\lambda}{\chi(\square)}$ goes along the same lines.

 By the definition of $\chi_i$, we have $\sum_{\square\in
\lambda}{\chi(\square)}=\sum_{i=1}^{\infty}{\frac{t_2^{i-1}-t_1\chi_i}{1-t_1}}$.
 Hence, we need to prove
$$[e_0,f_1]=-(1-t_1)^{-1}(1-t_2)^{-1}+\sum_{i=1}^{\infty}{\frac{t_2^{i-1}-t_1\chi_i}{1-t_1}}.$$
It is obvious for an empty diagram (and it is straightforward for
any $1$ row diagram $\lambda$).
 So it is enough to prove that
$[e_0,f_1]-\sum_{i=1}^{\infty}{\frac{t_2^{i-1}-t_1\chi_i}{1-t_1}}$
does not depend on $\lambda$.

 Choose $k\in \BN$, so that $\lambda_{k-1}=0$.
Then $\chi_i=t_1^{-1}t_2^{i-1}$ for $i\geq k-1$ and so~(\ref{KUKU})
implies:
$$\gamma_1=(1-t_1)^{-2}\left(
\sum_{i=1}^{k}{\frac{\chi_i(1-t_1t_2^{2-k}\chi_i)}{\chi_i-t_2^{k-1}}\prod_{1\leq
j\ne i}^{k-1}{\frac{(\chi_j-t_2\chi_i)(\chi_i-t_1t_2\chi_j)}
{(\chi_j-\chi_i)(\chi_i-t_1\chi_j)}}}\right.$$$$-\left.
\sum_{i=1}^{k}{\frac{\chi_i(1-t_1^2t_2^{2-k}\chi_i)}{\chi_i-t_1^{-1}t_2^{k-1}}\prod_{1\leq
j\ne i}^{k-1}{\frac{(\chi_i-t_2\chi_j)(\chi_j-t_1t_2\chi_i)}
{(\chi_i-\chi_j)(\chi_j-t_1\chi_i)}}}\right).$$

 The right-hand side is a rational expression in $\chi_i\ (1\leq i\leq k-2)$. The degree
of its numerator is not greater than that of the denominator plus 1,
while the only possible poles can happen at $\chi_i=\chi_j,\
\chi_i=t_1\chi_j$ (see the argument of (a)). These are simple poles
with zero residues. Hence, $\gamma_1$ is a linear function in
$\chi_i\ (1\leq i\leq k-2)$. It is straightforward to check that the
principal part of $\gamma_1$ equals
$\sum_{i=1}^{k-2}{-t_1(1-t_1)^{-1}\chi_i}$, and so
$[e_0,f_1]-\sum_{i=1}^{\infty}{\frac{t_2^{i-1}-t_1\chi_i}{1-t_1}}$
is constant. This completes the proof of
$[e_0,f_1]=-(1-t_1)^{-1}(1-t_2)^{-1}+\sum_{\square\in
\lambda}{\chi(\square)}$.
\end{proof}

\begin{corollary}
 Since $\gamma_0\equiv \frac{-1+t_1^{-1}t_2^{-1}}{(1-t_1)(1-t_2)(1-t_1^{-1}t_2^{-1})}$, we can define
$\phi^{+}_0:=-1,\ \phi^{-}_0:=-t_1^{-1}t_2^{-1}.$
\end{corollary}

 All the operators $\phi_s^{\pm}$ are
diagonalizable in the fixed point basis (according to
Proposition~\ref{diagonal}).


\noindent
 To prove~(\ref{4'}), it suffices to verify the following identity for any Young diagrams $\lambda,\
 \lambda'=\lambda+i_1$:
\begin{multline}\label{15}
(\phi^{+}_{i+3}e_j-\sigma_1\phi^{+}_{i+2}e_{j+1}+\sigma_2\phi^{+}_{i+1}e_{j+2}-\sigma_3\phi^{+}_ie_{j+3})_{[\lambda,\lambda{}']}\\=
(\sigma_3e_j\phi^{+}_{i+3}-\sigma_2e_{j+1}\phi^{+}_{i+2}+\sigma_1e_{j+2}\phi^{+}_{i+1}-e_{j+3}\phi^{+}_i)_{[\lambda,\lambda{}']}.
\end{multline}

 Since ${e_{j+k}}_{[\lambda,\lambda{}']}=(t_1^{\lambda_{i_1}}t_2^{i_1-1})^k{e_j}_{[\lambda,\lambda{}']}$ and $\phi_i$ is diagonalizable,
equality (\ref{15}) is equivalent to:
\begin{multline}
\label{6}
(\phi^{+}_{i+3}-\sigma_1\chi_1\phi^{+}_{i+2}+\sigma_2\chi_1^2\phi^{+}_{i+1}-\sigma_3\chi_1^3\phi^{+}_i)\mid_{\lambda{}'}\\=
(\sigma_3\phi^{+}_{i+3}-\sigma_2\chi_1\phi^{+}_{i+2}+\sigma_1\chi_1^2\phi^{+}_{i+1}-\chi_1^3\phi^{+}_i)\mid_{\lambda},
\end{multline}

\noindent where $j_1:=\lambda_{i_1}+1,\
\chi_1:=t_1^{j_1-1}t_2^{i_1-1},\ \phi^{+}_{<0}:=0$.
  First, we prove the analogous
equation for $\gamma_i$:
\begin{multline}
\label{6'}
(\gamma_{i+3}-\sigma_1\chi_1\gamma_{i+2}+\sigma_2\chi_1^2\gamma_{i+1}-\sigma_3\chi_1^3\gamma_i)\mid_{\lambda{}'}\\=
(\sigma_3\gamma_{i+3}-\sigma_2\chi_1\gamma_{i+2}+\sigma_1\chi_1^2\gamma_{i+1}-\chi_1^3\gamma_i)\mid_{\lambda}.
\end{multline}

\begin{proof}[Proof of equation~(\ref{6'})]
 Recalling the formula of $\gamma_i$, we consider every summand of $\gamma_i$ and compare its contribution to both sides of~(\ref{6'}).

\medskip

\noindent \textit{First case.} The summand of $\gamma_i$ corresponds
to the \textit{corner} $\square_{i_2,j_2}$ appearing in both sides
of~(\ref{6'}).

(a) Suppose $i_1<i_2,\ j_1>j_2$.

 Define $a:=j_1-j_2,\ b:=i_2-i_1,\ u:=t_1^a,\ v:=t_2^b,\ \chi_2:=t_1^{j_2-1}t_2^{i_2-1}$. So $\chi_1/\chi_2=uv^{-1}$.
 Then the contribution of $\square_{i_2,j_2}$ into
$(\gamma_{i+3}-\sigma_1\chi_1\gamma_{i+2}+\sigma_2\chi_1^2\gamma_{i+1}-\sigma_3\chi_1^3\gamma_i)\mid_{\lambda{}'}$
equals to
$$P_3\frac{(1-t_1^{-a}t_2^{b+1})(1-t_1^{a+1}t_2^{-b+1})}{(1-t_1^{-a}t_2^{b})(1-t_1^{a+1}t_2^{-b})}
\left(1-\sigma_1\left(\frac{\chi_1}{\chi_2}\right)+\sigma_2\left(\frac{\chi_1}{\chi_2}\right)^2-\sigma_3\left(\frac{\chi_1}{\chi_2}\right)^3\right)=$$$$
P_3(u-v)^{-1}(u-t_2v)(v-t_1t_2u)(v-t_2u)(v-t_1^{-1}t_2^{-1}u)v^{-3},$$
\medskip
while its contribution into
$(\sigma_3\gamma_{i+3}-\sigma_2\chi_1\gamma_{i+2}+\sigma_1\chi_1^2\gamma_{i+1}-\chi_1^3\gamma_i)\mid_{\lambda}$
equals to
$$P_3\frac{(1-t_1^{-a+1}t_2^{b+1})(1-t_1^{a}t_2^{-b+1})}{(1-t_1^{-a+1}t_2^{b})(1-t_1^{a}t_2^{-b})}
\left(\sigma_3-\sigma_2\left(\frac{\chi_1}{\chi_2}\right)+\sigma_1\left(\frac{\chi_1}{\chi_2}\right)^2-\left(\frac{\chi_1}{\chi_2}\right)^3\right)=$$$$
P_3(u-t_1t_2v)(u-v)^{-1}(v-t_2u)(u-t_2v)(u-t_1^{-1}t_2^{-1}v)v^{-3}.$$

 These two expressions coincide.

 (b) Case $i_1>i_2,\ j_1<j_2$ is completely analogous to (a).

\medskip

\noindent \textit{Second case.} The summand in the expression for
$\gamma_i$ corresponds to the \textit{hole} $\square_{i_2,j_2}$
which appears in both sides of~(\ref{6'}). In this case the argument
is analogous, since the expression for the summands of $\gamma$
corresponding to \textit{corners} and \textit{holes} differ only by
a sign.

\medskip
\noindent
 \textit{Third case.} Let us consider summands
occurring only in one side of~(\ref{6'}).
 Summands corresponding to deleting $\square_{i_1, j_1}$ in
the left-hand side of~(\ref{6'}) and to inserting
$\square_{i_1,j_1}$ in the right-hand side of~(\ref{6'}) are equal.
In the rest of the cases, both contributions are zero (we use the
argument that $t_1^{-1}, t_2^{-1}$ are roots of
$1-\sigma_1t+\sigma_2t^2-\sigma_3t^3$).
\end{proof}

 Now we are ready to verify equation~(\ref{6}).
\begin{proof} [Proof of equation~(\ref{6})]

 For $i>0$: ~(\ref{6}) is just~(\ref{6'}), proved above.

 The remaining (nontrivial) cases are $i=-3,\ -2,\ -1,\ 0$.
 According to~(\ref{6'}) and the relation between $\gamma_i$ and
$\phi^{\pm}_i$, we have to check only the following equalities:
$$\phi^{+}_0\mid_{\lambda{}'}=\phi^{+}_0\mid_{\lambda},\
\phi^{-}_0\mid_{\lambda{}'}=\phi^{-}_0\mid_{\lambda},\
\phi^{+}_1\mid_{\lambda{}'}=(\phi^{+}_1+(\sigma_1-\sigma_2)\chi_1
\phi^{+}_0)\mid_{\lambda},\
\phi^{-}_1\mid_{\lambda{}'}=(\phi^{-}_1+(\sigma_2-\sigma_1)\chi_1^{-1}
\phi^{-}_0)\mid_{\lambda}.$$

\noindent The first two are obvious, since $\phi^{\pm}_0$ are
constant.

\noindent Recalling
$\gamma_1\mid_{\mu}=-(1-t_1)^{-1}(1-t_2)^{-1}+\sum_{\square\in
\mu}{\chi(\square)}$ we get the third equality:
$$\phi_1^{+}\mid_{\lambda'}-\phi_1^{+}\mid_{\lambda}=(1-t_1)(1-t_2)(1-t_1^{-1}t_2^{-1})(\gamma_1\mid_{\lambda'}-\gamma_1\mid_{\lambda})=
(\sigma_1-\sigma_2)\chi_1 \phi^{+}_0.$$

\noindent The last one is proved in the same way.
\end{proof}

\medskip
  Equality ~(\ref{5'}) is proved in the same way as~(\ref{4'}).
Finally, we prove $\Phi^{+}(z)=\psi^{+}(z)$.

\medskip
 According to~(\ref{6}) we have
$$\Phi^{+}(z)(1-\sigma_1\chi_1z^{-1}+\sigma_2\chi_1^2z^{-2}-\sigma_3\chi_1^3z^{-3})\mid_{\lambda'}=
\Phi^{+}(z)(\sigma_3-\sigma_2\chi_1z^{-1}+\sigma_1\chi_1^2z^{-2}-\chi_1^3z^{-3})\mid_{\lambda}.$$

Thus
$$\Phi^{+}(z)\mid_{\lambda'}=\Phi^{+}(z)\mid_{\lambda}\cdot
\frac{(1-t_1^{-1}\chi_1z^{-1})(1-t_2^{-1}\chi_1z^{-1})(1-t_1t_2\chi_1z^{-1})}{(1-t_1\chi_1z^{-1})(1-t_2\chi_1z^{-1})(1-t_1^{-1}t_2^{-1}\chi_1z^{-1})}.$$

By induction $\Phi^{+}(z)\mid_{\lambda}=A\cdot \bc
(z)\mid_{\lambda}$, where $A$ is a coefficient of proportionality,
which is equal to $A=\Phi^{+}(z)\mid_{\emptyset}$ (here $\emptyset$
denotes an empty diagram). Explicitly:
$$A=\phi_0^{+}-(1-t_1^{-1}t_2^{-1})\sum_{i<0}{z^i}=-1-\frac{(1-t_1^{-1}t_2^{-1})z^{-1}}{1-z^{-1}}=
-\frac{1-t_1^{-1}t_2^{-1}z^{-1}}{1-z^{-1}}.$$

 So $\Phi^{+}(z)=\psi^{+}(z)$. Analogously one gets $\Phi^{-}(z)=\psi^{-}(z)$. $\square$

\medskip

\textrm{Theorem~\ref{main} is proved.} $\blacksquare$

\medskip

\section{Action of the shuffle algebra on $M$}
\label{section5}

In the previous section we constructed an action of the Ding-Iohara
algebra $A$ on $M$. Since parameters $q_1, q_2, q_3$ were not
generic (as we had $q_1q_2q_3=1$), Theorem~\ref{feigin1} does not
provide the representation of $S$ automatically. However, one can
get it just by writing the same formulas.

 Namely, let us define an action of $S$ on $M$ as follows.
For any $F\in S_n$ and any pair of Young diagrams $\lambda,\
\lambda'=\lambda+i_1+\cdots+i_n$ ($i_1\leq i_2\leq\cdots\leq i_n$),
we define the matrix coefficient
\begin{equation}
\label{shuffle-formula} F\mid_
{[\lambda,\lambda']}:=\frac{F(\chi_1,\ldots,\chi_n)}{\prod_{1\leq
a<b\leq
n}{\lambda(\chi_{i_a},\chi_{i_b})}}\prod_{k=1}^{n}e_{0[\lambda+i_1+\cdots+i_{k-1},\lambda+i_1+\cdots+i_k]},
\end{equation}

where $\chi_{i_k}$ is the character of the $k$th added box to
$\lambda$. All other matrix elements are zero.

\noindent The main result of this section is:

\begin{theorem}
\label{shuffle}

 Formula~(\ref{shuffle-formula}) provides a representation of the shuffle
 algebra $S$ on $M$.
\end{theorem}

 As a first step in the proof, we have the following statement:

 \begin{proposition}
 \label{order}

   If $\lambda, \lambda+j_1, \lambda+j_1+j_2,\ldots, \lambda'=\lambda+j_1+\cdots+j_n$ are Young diagrams, then
$$F\mid_{[\lambda,\lambda']}=\frac{F(\chi_1,\ldots,\chi_n)}{\prod_{1\leq a<b\leq
n}{\lambda(\chi_{j_a},\chi_{j_b})}}\prod_{k=1}^{n}e_{0[\lambda+j_1+\cdots+j_{k-1},\lambda+j_1+\cdots+j_k]},$$
where $\chi_{j_k}$ is the character of the $k$th added box to
$\lambda$ (first we add box $\#j_1$, then $\#j_2$, etc.). In other
words, condition $i_1\leq i_2\leq\cdots\leq i_n$ is irrelevant in
formula~(\ref{shuffle-formula}) for matrix coefficients.
 \end{proposition}

 \begin{proof}

  Since symmetric group is generated by transpositions, it suffices to check the statement only for them.
  However, this particular case follows from relation~(\ref{1}), in turn.
 \end{proof}

 \noindent{}Now we are ready to prove above theorem.
\begin{proof}[Proof of Theorem~\ref{shuffle}]

Let $F\in S_m,\ G \in S_n$, and let $\lambda,\
\lambda'=\lambda+j_1+\cdots+j_{m+n}$ be two Young diagrams. Then by
Proposition~\ref{order},
$$(F\circ G)\mid_{[\lambda, \lambda']}=\textrm{Sym}\left(\frac{G(\chi_1,\ldots,\chi_n)}{\prod_{1\leq
a<b\leq
n}{\lambda(\chi_{j_a},\chi_{j_b})}}\prod_{k=1}^{n}e_{0[\lambda+j_1+\cdots+j_{k-1},\lambda+j_1+\cdots+j_k]}\right.$$$$\times\left.\frac{F(\chi_{n+1},\ldots,\chi_{n+m})}{\prod_{n+1\leq
a<b\leq
n+m}{\lambda(\chi_{j_a},\chi_{j_b})}}\prod_{k=n+1}^{n+m}e_{0[\lambda+j_1+\cdots+j_{k-1},\lambda+j_1+\cdots+j_k]}\right)$$
$$=\textrm{Sym}\left(\frac{G(\chi_1,\ldots,\chi_n)F(\chi_{n+1},\ldots,\chi_{n+m})\prod_{k=1}^{n+m}e_{0[\lambda+j_1+\cdots+j_{k-1},\lambda+j_1+\cdots+j_k]}}{\prod_{1\leq
a<b\leq n}{\lambda(\chi_{j_a},\chi_{j_b})}\prod_{n+1\leq a<b\leq
n+m}{\lambda(\chi_{j_a},\chi_{j_b})}}\right).$$

\medskip

On the other hand: $$(G*F)(\chi_1,\ldots,\chi_{n+m})$$$$=
\textrm{Sym}\left(G(\chi_1,\ldots,\chi_n)F(\chi_{n+1},\ldots,\chi_{n+m})\prod_{1\leq
a\leq n< b\leq n+m}\lambda(\chi_{i_a},\chi_{i_b})\right).$$

 Thus applying Proposition~\ref{order}, we get
\footnote{\ Actually, we need $n!$ factor in
formula~(\ref{shuffle-formula}) if $\textrm{Sym}$ in
(\ref{star-product}) denotes the sum over all symmetric group
$\mathfrak{S}_{n+m}$, rather then over left cosets
$\mathfrak{S}_{n+m}/\mathfrak{S}_n\times \mathfrak{S}_m$.}
$$(F\circ G)\mid_{[\lambda, \lambda']}=(G*F)\mid_{[\lambda,
\lambda']}.$$

This completes the proof.
\end{proof}

 \noindent{}Let us now recall the collection of pairwise commuting elements $K_i\in S_i$ from Theorem~\ref{feigin2}:
$$K_2(x_1,x_2)=\frac{(x_1-q_1x_2)(x_2-q_1x_1)}{(x_1-x_2)^2},\ \ \ \ \ K_n(x_1,\ldots,x_n)=\prod_{1\leq i<j\leq
n}{K_2(x_i,x_j)}.$$
 Following result is crucial for the computations in the next section (recall that $q_1=t_1$):
\begin{lemma}
\label{Kn}

 For any $n$ indices $i_1<i_2<\cdots<i_n$, such that $\lambda':=\lambda+i_1+\cdots+i_n$ is a
Young diagram, we have:
$$K_n\mid_{[\lambda,\lambda']}=\prod_{1\leq a<b\leq
n}{\frac{(\chi_a-\chi_b)(\chi_b-t_1\chi_a)}{(\chi_a-t_2\chi_b)(\chi_a-t_1^{-1}t_2^{-1}\chi_b)}}\prod_{1\leq
r\leq
n}{e_{0[\lambda+i_1+\cdots+i_{r-1},\lambda+i_1+\cdots+i_r]}},$$
where $\chi_a=t_1^{\lambda_{i_a}}t_2^{{i_a}-1}.$
 All other matrix coefficients of $K_n$ are zero.
\end{lemma}

\emph{Remark:}
 We have constructed actions of Ding-Iohara and shuffle algebras
on $M$. The first one is purely geometric (given by operators $e_i,
f_i, \psi^{\pm}_i$), while the action of the shuffle algebra is
provided in algebraic terms. However, since elements $K_i$ belong to
the subalgebra generated by $S_1$ (Theorem 2.2), they are
geometrically represented since $S_1\subset \Xi(A_{+})$
(Theorem~\ref{feigin1}).

\section{Macdonald polynomials. Heisenberg algebra and vertex operators over it}
\label{section6}

\subsection{Macdonald polynomials}

In this subsection, we remind some basic facts about Macdonald
polynomials, following~\cite{m}.

 Recall that the algebra $\Lambda_F$ of symmetric functions over $F=\mathbb
Q(q,t)$ is freely generated by the power-sum symmetric functions
$\{p_k\}_{k\in \mathbb N}$; that is, $$\Lambda_F=F[p_1, p_2, \ldots\
].$$
 For any Young diagram $\lambda=(\lambda_1,\ldots,\lambda_k)=(1^{m_1}2^{m_2}\cdots)$ we define
$$p_{\lambda}:=p_{\lambda_1}\cdots p_{\lambda_k},\ \ \ \ \ z_{\lambda}:=\prod_{r\geq
1}{r^{m_r}m_r!}$$
  One introduces the Macdonald inner product $(\cdot,\cdot)_{q,t}$ on
  $\Lambda_F$ in the following way:
$$(p_{\lambda}, p_{\mu})_{q,t}=\delta_{\lambda,\mu}z_{\lambda}
\prod_{1\leq i\leq k}{\frac{1-q^{\lambda_i}}{1-t^{\lambda_i}}}.$$

\begin{definition}

 Macdonald polynomials $P_{\lambda}\in \Lambda_F$ are uniquely characterized by two
conditions:

 (a) $P_{\lambda}=m_{\lambda}+\textit{lower terms}$ (\textit{lower terms} stands for  $m_{\mu}$ with $\mu\prec \lambda$),

(b) $(P_{\lambda},P_{\mu})_{q,t}=0$ if $\lambda\ne \mu.$

\end{definition}

 \noindent{}Let $e_r$ be the $r$th elementary symmetric function. The following result is the \emph{Pieri formula}:

\begin{lemma}[\cite{m}, Section VI.6]

 We have the following multiplication rule: $$P_{\mu}e_r=\sum_{\lambda}\psi_{\lambda/\mu}P_{\lambda},$$ with the
sum taken over all $\lambda$ such that $\lambda/\mu$ is a vertical
$r$-strip and coefficients $\psi_{\lambda/\mu}$ equal to
\begin{equation}
\label{Pieri} \psi_{\lambda/\mu}=\prod
\frac{(1-q^{\mu_i-\mu_j}t^{j-i-1})(1-q^{\lambda_i-\lambda_j}t^{j-i+1})}{(1-q^{\mu_i-\mu_j}t^{j-i})(1-q^{\lambda_i-\lambda_j}t^{j-i})},
\end{equation}

where the product is taken over all pairs $(i, j)$ such that $i<j$
and $\lambda_i=\mu_i,\ \lambda_j=\mu_j+1$.

\end{lemma}

 \noindent{}In particular:
\begin{equation}
\label{Pieri2} \psi_{\mu+j/\mu}=
\prod_{i=1}^{j-1}\frac{(1-q^{\mu_i-\mu_j}t^{j-i-1})(1-q^{\mu_i-\mu_j-1}t^{j-i+1})}{(1-q^{\mu_i-\mu_j}t^{j-i})(1-q^{\mu_i-\mu_j-1}t^{j-i})}.
\end{equation}

\subsection{Fixed points via Macdonald polynomials}

 In this section we renormalize basis $[\lambda]$ of $M$ and operators $K_i\in S_i$, so that
 their action is consistent with an action of $e_i$ in the basis of
Macdonald polynomials. This is obtained via comparing the matrix
coefficients of $K_1$ and $e_1$.

 We define the normalized vectors $\langle \lambda
\rangle:=c_{\lambda}\cdot[\lambda]$, where
\begin{equation}
\label{renorm}
c_{\lambda}:=\left(-\frac{t_2}{1-t_2}\right)^{-|\lambda|}t_1^{\sum_{i}{\frac{\lambda_i(\lambda_i-1)}{2}}}\prod_{\square
\in
\lambda}{\left(1-t_1^{l(\square)}t_2^{-a(\square)-1}\right)}^{-1}.
\end{equation}
 First, we prove the following lemma:

\begin{lemma}
\label{normalization}

Let $\lambda=(\lambda_1,\ldots,\lambda_k)$. Then
$(1-t_1)(1-t_2)K_{1\langle \lambda, \lambda+j
\rangle}=\psi_{\lambda+j/\lambda}|_{q:=t_1, t:=t_2^{-1}}.$
\end{lemma}

\emph{Remark:} This lemma means that action of
$(1-t_1)(1-t_2)\Xi(e_0)$ in a renormalized basis is given by the
same formulas as an operator of multiplication by $e_1$ in the basis
of Macdonald polynomials. This condition defines a normalization
uniquely (up to a common factor).

\begin{proof}[Proof of Lemma~\ref{normalization}]
$\ $

 First, we have
$$\frac{c_{\lambda+j}}{c_{\lambda}}=-t_1^{\lambda_j}t_2^{-1}(1-t_2)
\prod_{\square \in
\Sigma_1(\square_{j,\lambda_j+1})}{\frac{1-t_1^{l(\square)}t_2^{-a(\square)-1}}{1-t_1^{l(\square)+1}t_2^{-a(\square)-1}}}\prod_{\square
\in
\Sigma_2(\square_{j,\lambda_j+1})}{\frac{1-t_1^{l(\square)}t_2^{-a(\square)-1}}{1-t_1^{l(\square)}t_2^{-a(\square)-2}}}
(1-t_2^{-1})^{-1}.$$

\medskip
 Let us rewrite each of these products in a more standard way:
$$\prod_{\square \in
\Sigma_2(\square_{j,\lambda_j+1})}{\frac{1-t_1^{l(\square)}t_2^{-a(\square)-1}}{1-t_1^{l(\square)}t_2^{-a(\square)-2}}}=
\prod_{i<j}{\frac{1-t_1^{\lambda_i-\lambda_j-1}t_2^{i-j}}{1-t_1^{\lambda_i-\lambda_j-1}t_2^{i-j-1}}},$$
$$\prod_{\square \in
\Sigma_1(\square_{j,\lambda_j+1})}{\frac{1-t_1^{l(\square)}t_2^{-a(\square)-1}}{1-t_1^{l(\square)+1}t_2^{-a(\square)-1}}}=
\frac{(1-t_2^{-1})}{(1-t_1^{\lambda_j}t_2^{j-k-1})}\prod_{i>j}^{k}{\frac{1-t_1^{\lambda_j-\lambda_i}t_2^{j-i-1}}{1-t_1^{\lambda_j-\lambda_i}t_2^{j-i}}}$$$$=
-(1-t_2^{-1})\prod_{i>j}^{\infty}{\frac{1-t_1^{\lambda_i-\lambda_j}t_2^{i-j+1}}{1-t_1^{\lambda_i-\lambda_j}t_2^{i-j}}}t_2t_1^{-\lambda_j}.
$$

\medskip
 Thus
$$\frac{c_{\lambda+j}}{c_{\lambda}}=(1-t_2)\prod_{i>j}^{\infty}{\frac{1-t_1^{\lambda_i-\lambda_j}t_2^{i-j+1}}{1-t_1^{\lambda_i-\lambda_j}t_2^{i-j}}}
\prod_{i<j}{\frac{1-t_1^{\lambda_i-\lambda_j-1}t_2^{i-j}}{1-t_1^{\lambda_i-\lambda_j-1}t_2^{i-j-1}}}.$$

 On the other hand, according to Proposition~\ref{eig'} we have
$$e_{0[\lambda,\lambda+j]}=(1-t_1)^{-1}\prod_{i\ne j}{\frac{1-t_1^{\lambda_i-\lambda_j}t_2^{i-j+1}}{1-t_1^{\lambda_i-\lambda_j}t_2^{i-j}}}.$$

 As for $\psi_{\lambda+j/\lambda}$, specializing formula~(\ref{Pieri2}) to $q:=t_1,\ t:=t_2^{-1}$, we get
$$\psi_{\lambda+j/\lambda}=
\prod_{i=1}^{j-1}\frac{(1-t_1^{\lambda_i-\lambda_j}t_2^{i-j+1})(1-t_1^{\lambda_i-\lambda_j-1}t_2^{i-j-1})}{(1-t_1^{\lambda_i-\lambda_j}t_2^{i-j})(1-t_1^{\lambda_i-\lambda_j-1}t_2^{i-j})}.$$

 Now it is straightforward to check that
 $$\psi_{\lambda+j/\lambda}=(1-t_1)(1-t_2)e_{0[\lambda,\lambda+j]}\cdot\frac{c_{\lambda}}{c_{\lambda+j}}=(1-t_1)(1-t_2)K_{1\langle \lambda, \lambda+j \rangle}.$$
\end{proof}

  Define constant factors $d_n:=\frac{(-t_1)^{n-1}}{(1-t_1)(1-t_2)}.$ Now we are ready to prove the main result:

\begin{theorem}\label {norm}

 For any Young diagrams $\mu\subset \lambda$, such that
$\lambda/\mu$ is a vertical $n$-strip with the boxes located in the
rows $j_1<\cdots<j_n$, we have: $\frac{1}{d_1\cdots
d_{n}}K_{n\langle \mu, \lambda
\rangle}=\psi_{\lambda/\mu}|_{q:=t_1,t:=t_2^{-1}}$.
\end{theorem}

\begin{proof}

 According to Lemma~\ref{Kn},
$$ K_{n[ \mu, \lambda ]}=\prod_{1\leq a<b\leq
n}{\frac{(\chi_a-\chi_b)(\chi_b-t_1\chi_a)}{(\chi_a-t_2\chi_b)(\chi_a-t_1^{-1}t_2^{-1}\chi_b)}}\prod_{1\leq
r\leq
n}{e_{0[\lambda-j_r-\cdots-j_n,\lambda-j_{r+1}-\cdots-j_{n}]}},$$
where $\chi_a=t_1^{\lambda_{j_a}-1}t_2^{{j_a}-1}.$
 Hence in the renormalized basis, $$K_{n\langle \mu, \lambda
\rangle}=\prod_{1\leq a<b\leq
n}{\frac{(\chi_a-\chi_b)(\chi_b-t_1\chi_a)}{(\chi_a-t_2\chi_b)(\chi_a-t_1^{-1}t_2^{-1}\chi_b)}}\prod_{1\leq
r\leq
n}{e_{0\langle\lambda-j_r-\cdots-j_n,\lambda-j_{r+1}-\cdots-j_n\rangle}}.$$
 On the other hand, specializing $q:=t_1,\ t:=t_2^{-1}$, formula~(\ref{Pieri}) reads as follows:
\begin{equation}
\label{Pieri3} \psi_{\lambda/\mu}=\prod
\frac{(1-t_1^{\mu_i-\mu_j}t_2^{i-j+1})(1-t_1^{\lambda_i-\lambda_j}t_2^{i-j-1})}{(1-t_1^{\mu_i-\mu_j}t_2^{i-j})(1-t_1^{\lambda_i-\lambda_j}t_2^{i-j})},
\end{equation}
where the product is over $\{(i,j)|i<j,\ i\notin
\{j_1,\ldots,j_n\}\ni j\}$. Explicitly (for $b$ fixed):

$$\prod_{j_1,\ldots,j_{b-1}\ne i<j_b}
\frac{(1-t_1^{\mu_i-\mu_{j_b}}t_2^{i-j_b+1})(1-t_1^{\lambda_i-\lambda_{j_b}}t_2^{i-j_b-1})}
{(1-t_1^{\mu_i-\mu_{j_b}}t_2^{i-j_b})(1-t_1^{\lambda_i-\lambda_{j_b}}t_2^{i-j_b})}\overset{\textit{Lemma
~\ref{normalization}}}=$$$$ \prod_{a<b}
{\frac{(1-t_1^{\lambda_{j_a}-\lambda_{j_b}+1}t_2^{j_a-j_b})(1-t_1^{\lambda_{j_a}-\lambda_{j_b}}t_2^{j_a-j_b})}
{(1-t_1^{\lambda_{j_a}-\lambda_{j_b}+1}t_2^{j_a-j_b+1})(1-t_1^{\lambda_{j_a}-\lambda_{j_b}}t_2^{j_a-j_b-1})}}\cdot
\frac{e_{0\langle \lambda-j_b-\cdots-j_n,
\lambda-j_{b+1}-\cdots-j_{n}\rangle}}{d_1}.$$

\medskip \noindent
Finally,
$$\prod_{a<b}
{\frac{(1-t_1^{\lambda_{j_a}-\lambda_{j_b}+1}t_2^{j_a-j_b})(1-t_1^{\lambda_{j_a}-\lambda_{j_b}}t_2^{j_a-j_b})}
{(1-t_1^{\lambda_{j_a}-\lambda_{j_b}+1}t_2^{j_a-j_b+1})(1-t_1^{\lambda_{j_a}-\lambda_{j_b}}t_2^{j_a-j_b-1})}}$$$$=
\prod_{a<b}{\frac{(\chi_b-t_1\chi_a)(\chi_b-\chi_a)}{(\chi_b-t_1t_2\chi_a)(\chi_b-t_2^{-1}\chi_a)}}=
\prod_{a<b}{\frac{(\chi_a-\chi_b)(\chi_b-t_1\chi_a)}{(\chi_a-t_2\chi_b)(\chi_a-t_1^{-1}t_2^{-1}\chi_b)}}(-t_1)^{1-b}.$$

This completes the proof of Theorem~\ref{norm} (since
$\frac{(-t_1)^{1-b}}{d_1}=\frac{1}{d_b}$).
\end{proof}

\subsection{Heisenberg action on $M$ via a vertex construction}
\label{Heisenberg1} $\ $

 Our results provide an isomorphism of $F$-vector spaces
$\Theta:M\iso \Lambda_F,\ \langle\lambda\rangle\mapsto P_\lambda$,
intertwining operators $\widetilde{K_i}:=\frac{1}{d_1\cdots
d_{i}}K_i\circlearrowright M$ with operators of multiplication by
$e_i\circlearrowright \Lambda_F$.

 Let us recall a beautiful identity of symmetric functions:
$$1+\sum_{i>0}{e_iz^i}=\exp\left(\sum_{i>0}{\frac{(-1)^{i-1}}{i}p_iz^i}\right).$$

 Thus identifying $M$ with a Fock space $\Lambda_F$, operators $\widetilde{K_i}$ can be viewed as \emph{vertex operators} over half of
the Heisenberg algebra $\{\mathfrak{h}_i\}_{i>0}$. In particular,
$\Theta$ intertwines $\mathfrak{h}_i$-action with operators of
multiplication by $p_i$ (for $i>0$). In other words, we get an
action of the positive part of the Heisenberg algebra. Analogous
considerations with operators $f_i$ instead of $e_i$ yield an action
of the negative part of the Heisenberg algebra. So in the end, we
get an action of the whole Heisenberg algebra on $M$.

\medskip
\emph{Remark:}
 Disadvantage of our approach is that we do not know explicit
formulas for $K_i$ in terms of
$x_1^{j_1}*x_1^{j_2}*\cdots*x_1^{j_i}$. \footnote{\ Such kind of
formulas were recently obtained by Negut (\cite{n}) for a whole
family of Heisenberg subalgebras.}

\subsection{Jack polynomials and equivariant cohomology as a specialization $q=t^{\alpha},\ t \to 1$.}

$\ \ \ \ $  Let us switch to considerations of equivariant
(localized) cohomology, instead of equivariant $K$-groups.
 Define $R:=\bigoplus_{n}\mathbb{H}^{2n}_{\mathbb T} (X^{[n]})\otimes_{\mathbb{H}_{\mathbb T}(\textrm{pt})}\on{Frac}(\mathbb{H}_{\mathbb T}(\textrm{pt}))$.
 According to the main result of~\cite{l}, there is a certain renormalization of the fixed point basis, such that
 isomorphism $\Delta: R\iso \Lambda_F$, defined by sending the renormalized basis of fixed
 points $\langle\lambda\rangle$ to Jack polynomials $J_{\lambda}^{(\alpha)}$, also intertwines operators $\{\mathfrak{h}_i\}_{i>0}$
 with operators of multiplication by $p_i$.

 On the other hand, it is known (see~\cite{m}) that Jack polynomials $J_{\lambda}^{(\alpha)}$ can be
obtained from the Macdonald polynomials $P_{\lambda}^{(q,t)}$ by
specializing $q:=t^{\alpha},\ t \to 1$. It is straightforward to
check that under this specialization, renormalizing
factors~(\ref{renorm}) become those of~\cite{l}, up to an additional
factor \footnote {\ This factor appears since our correspondences
slightly differ from those of Nakajima. While we use the whole
$P[\pm 1]$, Nakajima used only part of it, consisting of those
$(J_1,J_2)\in P[\pm 1]$, such that the quotient is supported at a
single point (which is automatic) with a zero $y$-coordinate (resp.,
$x$-coordinate).} (see formulas (2.12), (2.14) of~\cite{l} and note
that $l(\square),\ a(\square)$ are interchanged with our notations).
Since formulas in the fixed point basis of $R$ are just additive
analogues of the formulas for $M$, our approach provides the same
action of the Heisenberg algebra on $R$ as the classical one (via
higher correspondences $P[i]_{i\in \mathbb Z}$).

\section{Whittaker vector}
\label{section7}

 Let us consider the element $v=\sum_{n\geq 0}{\left[\mathcal{O}_{X^{[n]}}\right]}$ from a completion of $M$, denoted $M^\wedge$.
 We call it the \textit{Whittaker vector}, which is justified by the following theorem and~\cite{bf}, Section 2.30.

\begin{theorem}
\label{Whittaker}

 Consider $\widetilde{K}_{-n}:=\frac{1}{d_1\cdots
d_{n}}K_{-n}$ in the analogy with operators $\widetilde{K}_n$ from
Section~\ref{Heisenberg1}. Then for any $n\in \BN$ we have:
$\widetilde{K}_{-n}(v)=\widetilde{C}_n\cdot v$, where
$\widetilde{C}_n=\frac{(1-t_2)^n}{(1-t_2)(1-t_2^2)\cdots(1-t_2^n)}$.
\end{theorem}

 First we decompose $v$ in the fixed point basis.

\begin{proposition}
\label{technical}

 In the basis $[\lambda]$, vector $v$ is decomposed as follows:
$$v=\bigoplus_{\lambda}a_{\lambda}\cdot [\lambda], \ \
a_{\lambda}=\prod_{\square \in
\lambda}{\left(\left(1-t_1^{l(\square)+1}t_2^{-a(\square)}\right)\left(1-t_1^{-l(\square)}t_2^{a(\square)+1}\right)\right)^{-1}}.$$
\end{proposition}

\begin{proof}
 This is a straightforward consequence of the Bott-Lefschetz fixed point formula.
\end{proof}

\begin{proof}[Proof of Theorem~\ref{Whittaker}]

  We prove this theorem in two steps.

\noindent \textit{Step 1: Case $n=1$.}

  It suffices to prove $f_0(v)=K_{-1}(v)=C_1\cdot v$ with
$C_1=\frac{1}{(1-t_1)(1-t_2)}=d_1$, that is for any Young diagram
$\lambda$ the following identity holds: $C_1\cdot
a_{\lambda}=\sum_{j\leq k+1}{f_{0[\lambda+j,\lambda]}\cdot
a_{\lambda+j}}$, where $k$ denotes the length of $\lambda$.
According to Lemma~\ref{eigenv}, this is equivalent to
\begin{equation}
\label{nadoelo1} C_1=\sum_{j\leq
k+1}\frac{a_{\lambda+j}}{\chi_j\cdot a_{\lambda}}\prod_{s\in
\Sigma_1(\square_{j,\lambda_j+1})}{\frac{1-t_1^{l(s)+1}t_2^{-a(s)}}{1-t_1^{l(s)}t_2^{-a(s)}}}
\prod_{s\in
\Sigma_2(\square_{j,\lambda_j+1})}{\frac{1-t_1^{-l(s)}t_2^{a(s)+1}}{1-t_1^{-l(s)}t_2^{a(s)}}},
\end{equation}
 where $\chi_j=t_1^{\lambda_j}t_2^{j-1}$. Applying Proposition~\ref{technical}, we see that (\ref{nadoelo1}) is equivalent to
\begin{equation}
\label{nadoelo2} 1=\sum_{j\leq k+1}{\chi_j^{-1}}\prod_{s\in
\Sigma_1(\square_{j,\lambda_j+1})}{\frac{1-t_1^{-l(s)+1}t_2^{a(s)+1}}{1-t_1^{-l(s)}t_2^{a(s)+1}}}
\prod_{s\in
\Sigma_2(\square_{j,\lambda_j+1})}{\frac{1-t_1^{l(s)+1}t_2^{-a(s)+1}}{1-t_1^{l(s)+1}t_2^{-a(s)}}}.
\end{equation}

 As before, we rewrite both products from the right-hand side
 of~(\ref{nadoelo2}) in a uniform way:
$$\prod_{s\in
\Sigma_2(\square_{j,\lambda_j+1})}{\frac{1-t_1^{l(s)+1}t_2^{-a(s)+1}}{1-t_1^{l(s)+1}t_2^{-a(s)}}}=
\prod_{i<j}{\frac{1-t_1^{\lambda_i-\lambda_j}t_2^{1+i-j}}{1-t_1^{\lambda_i-\lambda_j}t_2^{i-j}}}=
\prod_{i< j}{\frac{\chi_j-t_2\chi_i}{\chi_j-\chi_i}},$$
$$\prod_{s\in
\Sigma_1(\square_{j,\lambda_j+1})}{\frac{1-t_1^{-l(s)+1}t_2^{a(s)+1}}{1-t_1^{-l(s)}t_2^{a(s)+1}}}=
\frac{1-t_2}{1-t_1^{-\lambda_j}t_2^{k-j+1}}\prod_{k\geq
i>j}{\frac{1-t_1^{\lambda_i-\lambda_j}t_2^{1+i-j}}{1-t_1^{\lambda_i-\lambda_j}t_2^{i-j}}}$$$$=
\prod_{k\geq i>j}{\frac{\chi_j-t_2\chi_i}{\chi_j-\chi_i}}\cdot
\frac{(1-t_2)\chi_j}{\chi_j-t_2^k}.$$

Hence~(\ref{nadoelo2}) is equivalent to the following identity:
\begin{equation}
\label{nadoelo3} 1=\sum_{j\leq k+1}{\prod_{k\geq i\ne
j}{\frac{\chi_j-t_2\chi_i}{\chi_j-\chi_i}}\cdot
\frac{1-t_2}{\chi_j-t_2^k}}.
\end{equation}

 Recalling the stabilizing condition $\chi_{k+1}=t_2^k$, this is
equivalent to
\begin{equation}
\label{nadoelo4} 1=\sum_{j\leq k+1}{\prod_{k+1\geq i\ne
j}{\frac{\chi_j-t_2\chi_i}{\chi_j-\chi_i}}\cdot
\frac{1-t_2}{\chi_j-t_2^{k+1}}}.
\end{equation}

 The right-hand side of~(\ref{nadoelo4}) is a rational function
$F(\chi_1,\ldots,\chi_{k+1})$, with
$\chi_j:=t_1^{\lambda_j}t_2^{j-1}$. In case
$\lambda_1=\cdots=\lambda_k=0$, that is $\chi_j=t_2^{j-1}$, we
obviously have $F(\chi_1,\ldots,\chi_{k+1})=1$. Since the degree of
the numerator of $F$ is not greater, then the degree of its
denominator, it suffices to check that $F$ does not have poles. This
verification is similar to the one from the proof of
Lemma~\ref{Heis1}. This finishes the proof of Step 1.

\medskip
\noindent \textit{Step 2: Case $n\geq 2.$}

Let $C_n=\frac{(-t_1)^{n(n-1)/2}}{(1-t_1)^n(1-t_2)\cdots(1-t_2^n)}$.
We will prove $K_{-n}(v)=C_n\cdot v$, which implies the result about
$\widetilde{K}_{-n}$, since $\widetilde{K}_{-n}=\frac{1}{d_1\cdots
d_{n}}K_{-n}$. The former equality is equivalent to (here $\lambda$
is arbitrary):
   $$C_n=\sum_{}{K_{-n[\lambda+i_1+\cdots+i_n,\lambda]}\cdot \frac{a_{\lambda+i_1+\cdots+i_n}}{a_{\lambda}}},$$
where the sum is over all sets of indices $i_1\leq \cdots\leq i_n$,
such that $\lambda+i_1+\cdots+i_n$ is a Young diagram.
 Similarly to Lemma~\ref{Kn}, we have
$$K_{-n[\lambda+i_1+\cdots+i_n,\lambda]}= \prod_{1\leq a<b\leq
n}{\frac{(\chi_a-\chi_b)(\chi_b-t_1\chi_a)}{(\chi_a-t_2\chi_b)(\chi_a-t_1^{-1}t_2^{-1}\chi_b)}}\prod_{1\leq
r\leq
n}{f_{0[\lambda+i_1+\cdots+i_{r},\lambda+i_1+\cdots+i_{r-1}]}},$$
where $\chi_a=t_1^{\lambda_{i_a}}t_2^{{i_a}-1}$. Hence
\begin{multline}
\label{nadoelo5} K_{-n[\lambda+i_1+\cdots+i_n,\lambda]}\cdot
\frac{a_{\lambda+i_1+\cdots+i_n}}{a_{\lambda}}=\prod_{1\leq
a<b\leq
n}{\frac{(\chi_a-\chi_b)(\chi_b-t_1\chi_a)}{(\chi_a-t_2\chi_b)(\chi_a-t_1^{-1}t_2^{-1}\chi_b)}}
\\ \times
\prod_{1\leq r\leq
n}{f_{0[\lambda+i_1+\cdots+i_{r},\lambda+i_1+\cdots+i_{r-1}]}\frac{a_{\lambda+i_1+\cdots+i_{r}}}
{a_{\lambda+i_1+\cdots+i_{r-1}}}}.
\end{multline}
 The latter terms were actually computed in \emph{Step 1}:
$$f_{0[\lambda+i_1+\cdots+i_{r},\lambda+i_1+\cdots+i_{r-1}]}\frac{a_{\lambda+i_1+\cdots+i_{r}}}
{a_{\lambda+i_1+\cdots+i_{r-1}}}=\frac{1}{(1-t_1)(1-t_2)}\prod_{k_r\geq
i\ne
i_r}{\frac{\chi_{i_r}-t_2\chi_i^{(r)}}{\chi_{i_r}-\chi_i^{(r)}}}\cdot
\frac{1-t_2}{\chi_{i_r}-t_2^{k_r}}=$$
$$\prod_{k+n\geq i\ne
i_r}{\frac{\chi_{i_r}-t_2\chi_i^{(r)}}{\chi_{i_r}-\chi_i^{(r)}}}\cdot
\frac{1}{(1-t_1)(\chi_{i_r}-t_2^{k+n})},$$ where
$\chi_i^{(r)}=t_1\chi_i$ for $i\in \{i_1,\ldots,i_{r-1}\}$, while
$\chi_i^{(r)}=\chi_i$ for $i\notin \{i_1,\ldots,i_{r-1}\}$, and
$k_r$ denotes the length of the diagram
$\lambda+i_1+\cdots+i_{r-1}$.

\noindent So the right-hand side of~(\ref{nadoelo5}) is equal to
$$\prod_{1\leq a<b\leq
n}{\frac{(\chi_a-\chi_b)(\chi_b-t_1\chi_a)}{(\chi_a-t_2\chi_b)(\chi_a-t_1^{-1}t_2^{-1}\chi_b)}}\prod_{1\leq
r\leq n}{\left(\frac{1}{(1-t_1)(\chi_{i_r}-t_2^{k+n})}
\prod_{k+n\geq i\ne
i_r}{\frac{\chi_{i_r}-t_2\chi_i^{(r)}}{\chi_{i_r}-\chi_i^{(r)}}}\right)}$$
$$={\prod_{1\leq r\leq n}{\left(\frac{(-t_1t_2)^{r-1}}{(1-t_1)(\chi_{i_r}-t_2^{k+n})}\prod_{k+n\geq j\ne
i_1,\ldots,i_n}{\frac{\chi_{i_r}-t_2\chi_j}{\chi_{i_r}-\chi_j}}\right)}}.$$

Define
\begin{multline}\label{G}
G(\chi_1,\ldots,\chi_{k+n}):= \sum_{i_1 < \cdots < i_n}{\prod_{1\leq
r\leq
n}{\left(\frac{(-t_1t_2)^{r-1}}{(1-t_1)(\chi_{i_r}-t^{k+n})}\prod_{k+n\geq
j\ne
i_1,\ldots,i_n}{\frac{\chi_{i_r}-t_2\chi_j}{\chi_{i_r}-\chi_j}}\right)}}.
\end{multline}

It is a rational function in variables $\chi_1,\ \ldots,\ \chi_k,\
\chi_{k+1}=t_2^k,\ \ldots,\ \chi_{k+n}=t_2^{k+n-1}$. The degree of
its numerator is not greater then the degree of the denominator.
Moreover, it has no poles (as before, it is crucial in the argument
that we have $\chi_{k+1}=t_2^k,\ldots,\ \chi_{k+n}=t_2^{k+n-1}$).

 Thus $G(\chi_1,\ldots,\chi_{k+n})$ is a constant function. Let us
calculate its value at $\{\chi_m=t_2^{m-1}\ \forall m\}$. All
summands of~(\ref{G}) vanish, except for the one corresponding to
$i_1=1,\ldots,\ i_n=n$. Hence
$$G=\frac{(-t_1t_2)^{n(n-1)/2}}{(1-t_1)^n(1-t_2^{k+n})\cdots(t_2^{n-1}-t_2^{k+n})}\prod_{n+1\leq j\leq
n+k}{\frac{(1-t_2^{j})\cdots(t_2^{n-1}-t_2^{j})}{(1-t_2^{j-1})\cdots(t_2^{n-1}-t_2^{j-1})}}=C_n.$$
\medskip
This completes a proof of the theorem.
\end{proof}

\begin{corollary}\label{uniq Wh 0}
 We have $f_{1}(v)=f_0(v)=C_1v$.
\end{corollary}

\begin{proof}
 Equality $f_1(v)=C_1v$ follows along the same lines as $f_0(v)=C_1v$ was proved.
\end{proof}

 Since the action of the negative half of the Heisenberg algebra
was constructed via the vertex operator construction
\begin{equation}
\label{vertex}
1+\sum_{i>0}{\widetilde{K}_{-i}z^i}=\exp\left(\sum_{i>0}{\frac{(-1)^{i-1}}{i}\mathfrak{h}_{-i}z^i}\right),
\end{equation}
we get the following corollary:

\begin{corollary}\label{Wh}
 We have $\mathfrak{h}_{-i}(v)=\alpha_i\cdot v$ for every $i\in \BN$, where $\alpha_i=\frac{(-1)^{i-1}(1-t_2)^i}{1-t_2^i}$.
\end{corollary}

\begin{proof}

 Firstly, Theorem~\ref{Whittaker} and~(\ref{vertex}) naturally imply that $v$ is a mutual eigenvector of all
$\{\mathfrak{h}_{-i}\}_{i\in \BN}$. Now we compute the eigenvalues
$\alpha_i$ explicitly.

  Using a combinatorial identity (a consequence of a $q$-binomial
formula)
\begin{equation}\label{q-binom1}
\sum_{i\geq 0}{\frac{z^i}{(1-t)(1-t^2)\cdots(1-t^i)}}=\prod_{i\geq
0} {\frac{1}{1-t^iz}}
\end{equation}
 and Theorem~\ref{Whittaker}, we get
$$\sum_{i>0}{\frac{(-1)^{i-1}}{i}\alpha_iz^i}=\ln\left(\sum_{i\geq
0}{\frac{(1-t_2)^i}{(1-t_2)(1-t_2^2)\cdots(1-t_2^i)}z^i}\right)$$$$=
-\ln(\prod_{i\geq 0}{(1-t_2^i(1-t_2)z)})=\sum_{i\geq 0}
{\sum_{j\geq 1}{\frac{t_2^{ij}(1-t_2)^jz^j}{j}}}=\sum_{j\geq 1}
{\frac{(1-t_2)^j}{1-t_2^j}\cdot \frac{z^j}{j}}.$$

Hence we have $\alpha_i=\frac{(-1)^{i-1}(1-t_2)^i}{1-t_2^i}$.
\end{proof}

 Let $[\mathfrak{h}_{-i},\mathfrak{h}_{i}]=\zeta_i\in \BC^*$, and define
$v_0:=\left[\mathcal{O}_{X^{[0]}}\right]\in M$. As an immediate
consequence of Corollary~\ref{Wh}, we get the following result:

\begin{proposition}

 We have $v=\exp(\sum_{i>0}{\frac{\alpha_i}{\zeta_i}\mathfrak{h}_i})v_0$.
\end{proposition}


\emph{Remark:}
 Completely analogously we can consider the \textit{Whittaker vector} $u$ in the completion of
 $R=\bigoplus_{n}\mathbb{H}^{2n}_{\mathbb T} (X^{[n]}) \otimes_{\mathbb{H}_{\mathbb T}(\mathrm{pt})}\on{Frac}(\mathbb{H}_{\mathbb{T}}(\mathrm{pt}))$.
 Similar arguments prove that $u=\exp(\frac{1}{\hbar\hbar'}\mathfrak{h}_1)u_0$,
 and hence $u$ is also an eigenvector with respect to the negative half of the Heisenberg algebra.

\subsection{Generalized Whittaker vector}\footnote {\ This section and the next one didn't appear in the published version.}
$\ $

 In this section we propose a slight generalization of the previous
 results. We start by introducing the following element $$w(x):=\sum_{n\geq
 0}{\left[{\Lambda^\bullet_x(\mathfrak{F})}_{X^{[n]}}\right]}\in M^\wedge[x].$$
 Specializing $x\in \BC$, it can be viewed as an element of $M^\wedge$.
Our main result is as follows:

\begin{theorem}\label{Whittaker generalized}

For every $n\in \BN$ we have

 (a) Operators $\widetilde{K}_{-n}$ act via $\widetilde{K}_{-n}(w(x))=\widetilde{C}_{x,n}\cdot w(x)$, where
     $\widetilde{C}_{x,n}=\frac{(1-t_2)^n\prod_{1\leq j\leq n}{(1-x t_2^{j-1})}}{\prod_{1\leq j\leq
     n}{(1-t_2^j)}}$,

 (b) Operators $\mathfrak{h}_{-n}$ act via
     $\mathfrak{h}_{-n}(w(x))=\alpha_{x,n}\cdot w(x)$, where
     $\alpha_{x,n}=\frac{(-1)^{n-1}(1-t_2)^n(1-x^n)}{1-t_2^n}$,
\medskip

 (c) We also have
$w(x)=\exp(\sum_{i>0}{\frac{\alpha_{x,i}}{\zeta_i}\mathfrak{h}_i})v_0$.
\end{theorem}

\begin{proof} [Sketch of the proof]$\ $

 (a) This statement is proved in the same way as
 Theorem~\ref{Whittaker}. In particular, formula~(\ref{nadoelo4}) from \emph{Step 1}, should be slightly changed to the following
 one:
 \begin{equation}
1-x=\sum_{j\leq k+1}{(1-x \chi_j)\prod_{k+1\geq i\ne
j}{\frac{\chi_j-t_2\chi_i}{\chi_j-\chi_i}}\cdot
\frac{1-t_2}{\chi_j-t_2^{k+1}}}.
\end{equation}

 This also yields a corresponding change of formulas in
 \emph{Step 2}.

 (b) This statement is proved the same way as Corollary~\ref{Wh} was deduced from Theorem~\ref{Whittaker}.
 We just need an updated version of identity~(\ref{q-binom1}), which now reads as follows
 \begin{equation}\label{q-binom2}
 1+\sum_{i\geq 1}{\left(z^i\prod_{1\leq j\leq i}{\frac{1-x t^{j-1}}{1- t^j}}\right)}=\prod_{i\geq 0} {\frac{1-xt^iz}{1-t^iz}}.
 \end{equation}
This is an easy consequence of (\ref{q-binom1}) and the formula
$\sum_{l=0}^k{\binom{k}{l}_qx^{k-l}(1-x)\cdots(1-q^{l-1}x)}=1$,
where  $k\in \NN$ and $\binom{k}{l}_q$ denotes a $q$-binomial
coefficient.

  (c) This is straightforward from (b).
\end{proof}

 \emph{Remark:} Note that specializing $x:=0$, we get $w(0)=v,\ \widetilde{C}_{0,n}=C_n,\ \alpha_{0,n}=\alpha_n$. Thus, Theorem~\ref{Whittaker
generalized} can be viewed as a straightforward generalization of
Theorem~\ref{Whittaker} and Corollary~\ref{Wh}.

  Specializing $x\in \BC$ we get the whole family of common eigenvectors of the commutative family of operators $\{K_{-n}\}_{n\in \BN}$.
 However, it turns out that Corollary~\ref{uniq Wh 0} determines $v$ uniquely:

\begin{proposition}
  If $v'\in M^\wedge$ satisfies $f_0(v')=f_1(v')=C_1v'$, then $v'=k\cdot v$ for some $k\in \BC$.
\end{proposition}

\begin{proof}
 Let $v'=\bigoplus_{\lambda}a'_{\lambda}\cdot [\lambda]$. It suffices to prove that fixing $a'_{\emptyset}$
 these conditions determine all $a'_{\lambda}$ uniquely.
 Our conditions can be rephrased as follows: for any Young diagram
 $\lambda$, we have
  $$C_1=\sum_{j\leq k+1}{f_{0[\lambda+j,\lambda]}\cdot \frac{a_{\lambda+j}}{a_\lambda}}=\sum_{j\leq k+1}{f_{1[\lambda+j,\lambda]}\cdot \frac{a_{\lambda+j}}{a_\lambda}},$$
  where $k$ denotes the length of $\lambda$. Since
  $f_{0[\lambda+j,\lambda]}=\chi_j^{-1}f_{1[\lambda+j,\lambda]}$
  with $\chi_j=t_1^{\lambda_j}t_2^{j-1}$ and all characters $\{\chi_j\}_{1\leq j\leq
  k+1}$ are pairwise distinct, an easy induction argument proves
  the statement (the induction is over $N:=|\lambda|$, while for a fixed $N$ it is over the dominance ordering of size $N$ Young diagrams).
   Note that existence of such $v'$ is guaranteed by Corollary~\ref{uniq Wh 0}.
\end{proof}

\subsection{Bilinear form}

 We finish this paper by determining a particular symmetric bilinear form on
 $M$, for which operators $e_i$ and $f_i$ are adjoint
 (compare this to Section 2.28 of~\cite{bf}).

\begin{definition}
 We define a symmetric $\BC(t_1,t_2)$-bilinear form $M\times M\overset{(\cdot,\cdot)}\longrightarrow \BC(t_1,t_2)$ as follows:

 $\bullet$ For $n\ne m$ we set $(M_n, M_m)=0$,

 $\bullet$ For any $n\in \BN_0$ and $\mathcal{G}_1, \mathcal{G}_2\in K^{\mathbb
 T}(X^{[n]})$, we define
 $(\mathcal{G}_1,\mathcal{G}_2):=[\mathrm{R\Gamma}(X^{[n]},\mathcal{G}_1\otimes \mathcal{G}_2\otimes \det(\mathfrak{F})^{-1})]$.
\end{definition}

 This definition is motivated by the following proposition:

\begin{proposition}\label{uniqueness}
 The above defined pairing satisfies the following two properties:

 (i) The normalizing condition: $(v_0,v_0)=1$.

 (ii) The adjoint condition: for any $i\in \BZ$ and $\mathcal{G}_1,\mathcal{G}_2\in M$ we have
     $(e_i(\mathcal{G}_1),\mathcal{G}_2)=(\mathcal{G}_1,f_i(\mathcal{G}_2))$.

 \noindent Conversely, these properties determine the pairing $(\cdot,\cdot)$ in a unique way.
\end{proposition}

\begin{proof}

 Statement (i) is obvious. Part (ii) follows from
 the projection formula together with an equality of line bundles $L\otimes q^*(\det(\mathfrak{F}))^{-1}=p^*(\det(\mathfrak{F}))^{-1}$ on the correspondence
 $P[1]$.

 The second part of the statement is proved by induction using equality $\sum_{i\in \BZ}{e_i(M_n)}=M_{n+1}$. In particular,
 we automatically get $(M_n,M_m)=0$ for $n\ne m$.
\end{proof}

 Let us now compute this form $(\cdot,\cdot)$ explicitly in the fixed point basis of $M$.

\begin{proposition}
In the fixed point basis $\{[\lambda]\}$ of $M$ we have:

(i) $([\lambda],[\lambda'])=0$ for $\lambda\ne \lambda'$,

(ii)
$([\lambda],[\lambda])=a_{\lambda}^{-1}t_1^{-\sum_i{\frac{\lambda_i(\lambda_i-1)}{2}}}t_2^{-\sum_i{(i-1)\lambda_i}}=
a_{\lambda}^{-1}t_1^{-\sum_i{\frac{\lambda_i(\lambda_i-1)}{2}}}t_2^{-\sum_j{\frac{\check{\lambda}_j(\check{\lambda}_j-1)}{2}}}$,
where $a_\lambda$ were defined in Proposition~\ref{technical} and
$\check{\lambda}$ is the transpose of the Young diagram $\lambda$.
\end{proposition}

\begin{proof}
 It is straightforward to check that formulas (i--ii) define a
symmetric $\BC(t_1,t_2)$-bilinear form on $M$, satisfying conditions
(i--ii) of Proposition~\ref{uniqueness} (computations are similar to
those in Section 6). Aforementioned proposition implies the result.
\end{proof}

 We finish our discussion of $(\cdot,\cdot)$ by comparing it to the Macdonald inner product $(\cdot,\cdot)_{q,t}$ from Section 6.
 Recalling the renormalized basis $\langle\lambda\rangle=c_\lambda[\lambda]$ of $M$, we have the following result:
\begin{proposition}
 For any diagram $\lambda$ we have
 $(\langle\lambda\rangle,\langle\lambda\rangle)=(P_\lambda,P_\lambda)_{q:=t_1,t:=t_2^{-1}}\cdot(\frac{-(1-t_2)^2}{t_2})^{|\lambda|}$.
\end{proposition}
\begin{proof}
  According to~\cite{m}: $(P_\lambda,P_\lambda)_{q,t}=\prod_{\square\in \lambda}{\frac{1-q^{l(\square)+1}t^{a(\square)}}{1-q^{l(\square)}t^{a(\square)+1}}}$.
  The rest is straightforward.
\end{proof}

 In other words, isomorphism $M\iso \Lambda_F$ from Section 6,
 sending $\langle\lambda\rangle$ to
 $P_\lambda$, intertwines a bilinear form $(\cdot,\cdot)_{q,t}$ on $\Lambda_F$ and
 a bilinear form $(\cdot,\cdot)'$ on $M$, where the latter one is characterized uniquely by $(v_0,v_0)'=1$
 and $(e_i(\mathcal{G}_1),\mathcal{G}_2)'=\frac{-t_2}{(1-t_2)^2}\cdot(\mathcal{G}_1,f_i(\mathcal{G}_2))$
 for any $\mathcal{G}_1,\mathcal{G}_2\in M$ and $i\in \BZ$.

\end{document}